\documentclass{article}

\usepackage{natbib}
\usepackage[T1]{fontenc}    
\usepackage{hyperref}       
\usepackage{url}            
\usepackage{booktabs}       
\usepackage{amsfonts}       
\usepackage{nicefrac}       
\usepackage{microtype}      
\usepackage{xcolor}         
\usepackage{amsmath}
\usepackage{amssymb}
\usepackage{bm} 
\usepackage{float}
\usepackage{geometry}
\usepackage{amsthm}
\usepackage{mathrsfs}
\usepackage{enumerate}
\usepackage{cleveref}
\usepackage{algorithm}
\usepackage{algorithmic}
\usepackage{graphicx}  
\usepackage{subfigure}

\title{A More Stable Accelerated Gradient Method Inspired by Continuous-Time Perspective}
\date{}

\author{Yasong Feng\footnote{\texttt{ysfeng20@fudan.edu.cn}}\qquad
Weiguo Gao\footnote{\texttt{wggao@fudan.edu.cn}}\qquad
}

\newcommand{\df}{\rm{d}}
\newtheorem{lemma}{Lemma}

\newtheorem{rmk}{Remark}
\newtheorem{theo}{Theorem}

\newtheorem{prop}{Proposition}

\numberwithin{equation}{section}

\newcommand{\e}{\mathbb{E}}

\newcommand{\tr}{\mathsf{T}}

\begin{document}

\maketitle

\begin{abstract}
  Nesterov's accelerated gradient method (NAG) is widely used in problems with machine learning background including deep learning, and is corresponding to a continuous-time differential equation. From this connection, the property of the differential equation and its numerical approximation can be investigated to improve the accelerated gradient method. In this work we present a new improvement of NAG in terms of stability inspired by numerical analysis. We give the precise order of NAG as a numerical approximation of its continuous-time limit and then present a new method with higher order. We show theoretically that our new method is more stable than NAG for large step size. Experiments of matrix completion and handwriting digit recognition demonstrate that the stability of our new method is better. Furthermore, better stability leads to higher computational speed in experiments.
\end{abstract}

\section{Introduction}
Optimization is a core component of statistic and machine learning problems. Recently, gradient-based algorithms are widely used in such optimization problems due to its simplicity and efficiency for large-scale situations. For solving convex optimization problem
\begin{align*}
\min_{x\in\mathbb{R}^d}F(x),
\end{align*}
where $F(x)$ is convex with Lipschitz gradient, it is known that the objective function converges at a rate of $\mathcal{O}(n^{-1})$ when applying vanilla gradient descent.

Alternatively, \citet{nes} provided a more efficient first-order algorithm than gradient descent, which may take the following form: starting with $x_0=x_1$,
\begin{align}\label{nes}
\begin{split}
&y_n=x_n+\frac{n-3}{n}(x_n-x_{n-1}),\\
&x_{n+1}=y_{n}-s\nabla F(y_{n})
\end{split}
\end{align}
for $n\geq1$. It is shown that under abovementioned conditions, Nesterov's accelerated gradient method (NAG) converges at the optimal rate $\mathcal{O}(n^{-2})$.

Accelerated gradient method has been successful in training deep and recurrent neural networks \citep{13} and is widely used in problems with machine learning background to avoid sophisticated second-order methods \citep{09,ji,11}. To provide more theoretical understanding, an important research topic of NAG is to find an explanation of the acceleration. \citet{su} took the step size $s\to0$ and showed that the continuous-time limit of NAG is the following ODE
\begin{equation}\label{ode}
\ddot{x}+\frac{3}{t}\dot{x}+\nabla F(x)=0.
\end{equation}
This differential equation was used as a tool for analyzing and generalizing Nesterov's scheme. 

From this continuous-time perspective, progresses have been made in presenting new continuous frameworks, for example, see \citet{wibi,lyap,highr} and \citet{dis}. These new frameworks give more insight on the acceleration phenomenon. Furthermore, discretization methods are provided to derive optimization algorithms from a certain continuous-time framework. \cite{zhang} showed that discrete the ODE (\ref{ode}) using a Runge-Kutta integrator can achieve acceleration. \cite{sdhr} showed that applying symplectic scheme to a ``high-resolution" version of limit ODE can also achieve acceleration. \citet{rgd} presented a new method by discretizing a ``dissipative relativistic system" that can be viewed as a generalization of NAG and heavy ball method.

Despite the successful use in practice, Nesterov's accelerated gradient method is relatively poor in stability, that is, it requires a small step size to ensure convergence when applied to an objective function with large gradients \citep{aggr,rgd}, which limits the computational speed. To improve the stability, we present a new method motivated by numerical analysis. We prove that NAG is actually a numerical approximation of ordinary differential equation (\ref{ode}), so the properties in numerical analysis can be adopted to improve it. The order of a numerical approximation method is closely related with its stability. For example, the absolutely stable region of Runge-Kutta method is enlarged when the order is increased from one to four \citep{rkab}. In this work we give the precise order of Nesterov's accelerated gradient method as a numerical approximation of the limit ODE. Then inspired by the abovementioned property of Runge-Kutta method, we present a new higher order method (Algorithm \ref{sag}), which we call \emph{stabilized accelerated gradient} (SAG). Compared with NAG, we use one more step to calculate the starting point $\bm{Y}_k$ for every iteration, so our SAG approximates the ordinary differential equation (\ref{ode}) better. We emphasis here that in many cases (including the two experiments in our work), the most computationally expensive step in gradient-based algorithms is the calculation of the gradients, so our modification will not cause inefficiency. Moreover, using absolute stability theory, we prove that our SAG is more stable than Nesterov's method.
\begin{algorithm}
	\caption{Stabilized Accelerated Gradient (SAG)}
	\label{sag}
	\begin{algorithmic}
		\STATE \textbf{Input: }step size $s$.
		\STATE \textbf{Initial value: }$\bm{X}_2=\bm{X}_1=\bm{X}_0$.
		\STATE \textbf{$\bm{(k-1)}$th iteration $(k\geq2)$.} Compute
		\begin{align*}
		&\bm{Y}_k=\frac{10k^2+9k+6}{4k^2+8k}\bm{X}_k-\frac{4k^2+3}{2k^2+4k}\bm{X}_{k-1}+\frac{2k-1}{4k+8}\bm{X}_{k-2},\\
		&\bm{Z}_k=\frac{2k-3}{k}\bm{X}_k-\frac{k-3}{k}\bm{X}_{k-1},\\
		&\bm{X}_{k+1}=\bm{Y}_k-\frac{ks}{2k+4}\nabla F(\bm{Z}_k).
		\end{align*}
	\end{algorithmic}
\end{algorithm} 

Based on abovementioned theoretical results, we try to take advantage of the new method in more practical problems. We apply SAG to a large-scale matrix completion problem. We combine SAG with proximal operator \citep{prox} into a new algorithm, which we call SFISTA. We show that the performance of our SFISTA is better than FISTA \citep{fista} and accelerated proximal gradient method \citep{prox} because of its better stability.

Furthermore, we show the advantage of our SAG by a handwriting digit recognition task. We train a convolutional neural network using SAG and NAG on \emph{MNIST} \citep{mni} dataset. This experiment also verifies that our SAG is more stable with large step size.

This paper is organized as follows. In section $2$,  We prove that NAG is a numerical approximation of ordinary differential equation (\ref{ode}) and give the precise order. In section $3$, we present our higher order method SAG and prove its better stability. In section $4$ and $5$, we apply our SAG to matrix completion and handwriting digit recognition tasks.

\section{The precise order of Nesterov's method as a numerical approximation of its continuous-time limit}
We refer to $x(t)$ as the solution of differential equation (\ref{ode}). \citet{su} proved the existence and uniqueness of such a solution. In this section, We give the precise order of NAG as a numerical approximation of $x(t)$.

We substitute the first equation in Nesterov's scheme (\ref{nes}) to the second one and substitute the relation $s=h^2$ to get
\begin{equation}\label{ite1}
x_{n+1}=x_n+\frac{n-3}{n}(x_n-x_{n-1})-h^2\cdot \nabla F\left(x_n+\frac{n-3}{n}(x_n-x_{n-1})\right).
\end{equation}
Inspired by the ansatz $x_n\approx x(n\sqrt{s})=x(nh)$ \citep{su}, we consider the convergence between $x_n$ and $x(nh)$. More precisely, we show that for fixed time $t$, $x_n$ converges to $x(t)$ with order $\mathcal{O}(h\ln\frac{1}{h})$, where $n=\frac{t}{h}$.
\subsection{Truncation error}
Firstly, we consider the following ``truncation error":
\begin{align}\label{rn}
\begin{split}
L[x(t);h]=&x(t+h)-\frac{2t-3h}{t}x(t)+\frac{t-3h}{t}x(t-h)+\\&h^2\nabla F\left(x(t)+\frac{t-3h}{t}\left(x(t)-x(t-h)\right)\right).
\end{split}
\end{align}
(\ref{rn}) is obtained from (\ref{ite1}) by replacing $x_{n+1},\;x_n,\;x_{n-1}$ with $x(t+h),\;x(t),\;x(t-h)$ and substituting the relation $n=\frac{t}{h}$. Our first result is the order of truncation error $L[x(t);h]$.

\begin{theo}\label{trun}
	Assume $\nabla F$ satisfies $L$-Lipschitz condition, and solution $x(t)$ of the derived differential equation (\ref{ode}) has a continuous third derivative. For fixed time $t$, the truncation error (\ref{rn}) satisfies
	\begin{equation*}
	L[x(t);h]=\mathcal{O}(h^3).
	\end{equation*}
\end{theo}
\subsection{Approximation theorem}
Theorem \ref{trun} coincides with derivation of ODE (\ref{ode}) \citep{su} and shows the size of error caused by a single iteration when the starting point is just on $x(t)$. Then we add up these errors to give the precise order of $x_n$ converging to its limit dynamic system.
\begin{theo}\label{conv}
	Under conditions in Theorem \ref{trun}, for fixed time $t$ and $n=t/h$, $x_n$ converges to $x(t)$ at a rate of $\mathcal{O}(h\ln\frac{1}{h})$ if $x_0=x(0)$ and $x_1=x(h)$.
\end{theo}
For the proof of Theorem \ref{conv}, we need the following two lemmas.
\begin{lemma}\citep{gron}\label{gron}
	For constant $\alpha,\;\beta>0$ and positive sequence $\{\eta_n\}_{n\geq0}$ satisfying
	\[\eta_n\leq\beta+\alpha\sum_{i=0}^{n-1}\eta_i,\quad\forall n>0,\]
	the following inequality holds
	\[\eta_n\leq e^{\alpha n}(\beta+\alpha\eta_0).\]
\end{lemma}
The above lemma is a classic result and referred to as discrete Gronwall inequality.
\begin{lemma}\label{mai}
	We define matrices $\bm{C}_n$ and $\bm{D}_{n,l}$ as
	\begin{align*}
	&\bm{C}_n=\begin{pmatrix}\frac{2n-1}{n+1}&&-\frac{n-2}{n+1}\\1&&0\end{pmatrix},\\
	&\bm{D}_{n,l}=\bm{C}_n\bm{C}_{n-1}\cdots\bm{C}_{n-l+1},
	\end{align*}
	where $n\geq0$ and $0<l\leq n+1$. In addition, we set $\bm{D}_{n,0}=\bm{I}_2$. Then there exist positive constants $M,\;M_3$ such that for all $n$, the following two inequalities hold, where the matrix norm is 2-norm.
	\begin{align*}
	\sup_{0\leq l\leq n+1}\|\bm{D}_{n,l}\|&\leq Mn,\\
	\bm{D}_{n,n+1}&\leq M_3.
	\end{align*}
\end{lemma}
The proof of Theorem \ref{trun}, Lemma \ref{mai} and Theorem \ref{conv} are presented in Appendix \ref{a1}, \ref{a2} and \ref{a3}.

We notice that \cite{su} also presented an approximation result (Theorem 2 in \citet{su}). However, our new approach is based on the order of truncation error and follows the classical path in numerical analysis that from consistence to convergence. This approach allows us to present the precise order of NAG as an approximation of the limit ODE, which is new to our knowledge. Based on this new result, we present our higher order method.

\section{Stabilized accelerated gradient method}
From Theorem \ref{conv}, NAG can be viewed as a numerical approximation of the limit ODE (\ref{ode}). Therefore, the property of this limit differential equation and its numerical approximation can be investigated to improve the accelerated gradient method. In this section, we present a new way to improve the stability of NAG. It is well known that the absolutely stable region of Runge-Kutta method is enlarged when the order increases from one to four. Inspired by this fact, we present our new higher order method SAG. Here we show the derivation of SAG and then prove that our SAG is more stable than NAG.
\subsection{Derivation of the new method and analysis of truncation error} \label{nmtd}
We consider a recurrence relation with the following form
\begin{align*}
\begin{split}
\sum_{i=1}^4\left(\alpha_i+\frac{\beta_i}{n}+\frac{\gamma_i}{n^2}\right)x_{n+2-i}=-h^2\nabla F\left(x_n+\frac{n-3}{n}(x_n-x_{n-1})\right),
\end{split}
\end{align*}
where $\{\alpha_i\}$, $\{\beta_i\}$ and $\{\gamma_i\}$ are to be determined. For this scheme, the truncation error becomes
\begin{align*}
\begin{split}
L[x(t);h]=&\sum_{i=1}^4\left(\alpha_i+\frac{\beta_ih}{t}+\frac{\gamma_ih^2}{t^2}\right)x(t+(2-i)h)\\&+h^2\nabla F\left(x(t)+\frac{t-3h}{t}(x(t)-x(t-h))\right).
\end{split}
\end{align*}
For the gradient term, we expand $x(t-h)$ to first order to get
\begin{align*}
\nabla F\left(x(t)+\frac{t-3h}{t}(x(t)-x(t-h))\right)=&-hx^{(3)}(t)-\left(\frac{3h}{t}+1\right)x^{(2)}(t)\\
&+\left(\frac{3h}{t^2}-\frac{3}{t}\right)x^{(1)}(t)+\mathcal{O}(h^2).
\end{align*}
Then we substitute this expansion for the gradient term in the truncation error. Calculation shows that the terms with order less than four of the truncation error will be eliminated if and only if the following equations are satisfied,
\[\left\{\begin{aligned}\alpha_1&=2\\\alpha_2&=-5\\\alpha_3&=4\\\alpha_4&=-1\end{aligned}\right.,\qquad\left\{\begin{aligned}\beta_1&=\frac{9}{2}-k\\\beta_2&=-6+3k\\\beta_3&=\frac{3}{2}-3k\\\beta_4&=k\end{aligned}\right.,\qquad\left\{\begin{aligned}\gamma_1&=m_1\\\gamma_2&=-\frac{3m_1+m_2+3}{2}\\\gamma_3&=m_2\\\gamma_4&=\frac{m_1-m_2+3}{2}\end{aligned}\right.,\]
where $k,\;m_1,\;m_2$ can be chosen randomly. Then we set $k=\frac{1}{2}$, $m_1=0$ and $m_2=3$ to get our \emph{stabilized accelerated gradient method} (Algorithm \ref{sag})
\begin{align}\label{ite2e}
\begin{split}
x_{n+1}=&\frac{10n^2+9n+6}{4n^2+8n}x_n-\frac{4n^2+3}{2n^2+4n}x_{n-1}+\frac{2n-1}{4n+8}x_{n-2}\\&-\frac{n}{2n+4}h^2\nabla F\left(\frac{2n-3}{n}x_n-\frac{n-3}{n}x_{n-1}\right).
\end{split}
\end{align}
For truncation order of SAG, we have the following theorem. The abovementioned procedure is presented in detail in Appendix \ref{a4}, as proof of Theorem \ref{t_rnorder}.
\begin{theo}\label{t_rnorder}
	If $\nabla F$ has continuous second order derivative, the first and second derivative are bounded, and $x(t)$ has continuous fourth derivative, then for fixed $t$, truncation error of SAG (\ref{ite2e}) satisfies 
	\begin{equation*}
	L[x(t);h]=\mathcal{O}(h^4).
	\end{equation*}
\end{theo}
\subsection{Better stability of SAG}
In this subsection, we verify the better stability of the new method. The stability here refers to the feasible region of step size. It is shown that step size should be chosen in absolute stable region to ensure stable performance of the algorithm \citep{abs}. Therefore, a more stable method should have a larger absolute stable region, so that its effectiveness is ensured when the step size is large. 

We emphasize here that since absolute stable theory is corresponding to accumulated error, it can be widely used in iterative algorithms, although this concept was orginally proposed for numerical integrators. Similar analyses are applied to different optimization methods in \citet{su} and \citet{rgd}.

Firstly, recall the scheme of SAG
\begin{align*}
x_{n+1}=&\frac{10n^2+9n+6}{4n^2+8n}x_n-\frac{4n^2+3}{2n^2+4n}x_{n-1}+\frac{2n-1}{4n+8}x_{n-2}\\&-\frac{n}{2n+4}s\nabla F\left(\frac{2n-3}{n}x_n-\frac{n-3}{n}x_{n-1}\right).
\end{align*}
We use notations $\bm{X}_{n+1}=(x_{n+1},x_n,x_{n-1})^\tr$ and $\bm{X}_{n}=(x_{n},x_{n-1},x_{n-2})^\tr$, then the scheme can be written as
\[\bm{X}_{n+1}=F^{(n)}(\bm{X}_n)=\begin{pmatrix}F_1^{(n)}(\bm{X}_n)\\F_2^{(n)}(\bm{X}_n)\\F_3^{(n)}(\bm{X}_n)\end{pmatrix},\]
where $F_1^{(n)}(x,y,z)=\frac{10n^2+9n+6}{4n^2+8n}x-\frac{4n^2+3}{2n^2+4n}y+\frac{2n-1}{4n+8}z-\frac{n}{2n+4}s\nabla F\left(\frac{2n-3}{n}x-\frac{n-3}{n}y\right)$, $F_2^{(n)}(x,y,z)=x$ and $F_3^{(n)}(x,y,z)=y$.

In a real SAG run, we calculate a numerical approximation $\bar{\bm{X}_n}$ of each $\bm{X}_n$ in turn. Assume that there is an error $\bm{e}=(e_1,e_2,e_3)^\tr$ between $\bar{\bm{X}_n}$ and $\bm{X}_n$. Then in the next iteration, we have 
\begin{align}\label{erreq}
    \bar{\bm{X}}_{n+1}=F^{(n)}(\bar{\bm{X}_n})=\begin{pmatrix}F_1^{(n)}(\bm{X}_n+\bm{e})\\F_2^{(n)}(\bm{X}_n+\bm{e})\\F_3^{(n)}(\bm{X}_n+\bm{e})\end{pmatrix}.
\end{align}
For the first entry of RHS, we have
\begin{align*}
    F^{(n)}_1(\bm{X}_n+\bm{e})=&\frac{10n^2+9n+6}{4n^2+8n}(x_n+e_1)-\frac{4n^2+3}{2n^2+4n}(x_{n-1}+e_2)+\frac{2n-1}{4n+8}(x_{n-2}+e_3)\\
    &-\frac{n}{2n+4}s\nabla F\left(\frac{2n-3}{n}(x_n+e_1)-\frac{n-3}{n}(x_{n-1}+e_2)\right)\\
    \approx&\frac{10n^2+9n+6}{4n^2+8n}(x_n+e_1)-\frac{4n^2+3}{2n^2+4n}(x_{n-1}+e_2)+\frac{2n-1}{4n+8}(x_{n-2}+e_3)\\
    &-\frac{n}{2n+4}s\nabla F\left(\frac{2n-3}{n}x_n-\frac{n-3}{n}x_{n-1}\right)\\
    &-\frac{n}{2n+4}s\nabla^2 F\left(\frac{2n-3}{n}x_n-\frac{n-3}{n}x_{n-1}\right)\cdot\left(\frac{2n-3}{n}e_1-\frac{n-3}{n}e_2\right)\\
    =&F^{(n)}_1(\bm{X}_n)+\left(\frac{10n^2+9n+6}{4n^2+8n}-\frac{2n-3}{2n+4}s\nabla^2 F\right)e_1-\left(\frac{4n^2+3}{2n^2+4n}-\frac{n-3}{2n+4}s\nabla^2 F\right)e_2\\&+\frac{2n-1}{4n+8}e_3,
\end{align*}
where in the second line we use Taylor approximation, and in the third line we use $\nabla^2F$ to denote $\nabla^2 F\left(\frac{2n-3}{n}x_n-\frac{n-3}{n}x_{n-1})\right)$. 

Then (\ref{erreq}) implies that
\begin{align*}
    \bar{\bm{X}}_{n+1}&\approx \bm{X}_{n+1}+\begin{pmatrix}\left(\frac{10n^2+9n+6}{4n^2+8n}-\frac{2n-3}{2n+4}s\nabla^2 F\right) && -\left(\frac{4n^2+3}{2n^2+4n}-\frac{n-3}{2n+4}s\nabla^2 F\right) && \frac{2n-1}{4n+8}\\
    1 && 0 && 0\\
    0 && 1 && 0\end{pmatrix}
    \begin{pmatrix}
    e_1\\e_2\\e_3
    \end{pmatrix}\\
    &\triangleq\bm{X}_{n+1}+\bm{P}\bm{e}.
\end{align*}

The characteristic equation of matrix $\bm{P}$ is \[\lambda^3-\left(\frac{10n^2+9n+6}{4n^2+8n}-\frac{2n-3}{2n+4}s\nabla^2 F\right)\lambda^2+\left(\frac{4n^2+3}{2n^2+4n}-\frac{n-3}{2n+4}s\nabla^2 F\right)\lambda-\frac{2n-1}{4n+8}=0.\]
For large $n$, we ignore the high order terms and the characteristic equation becomes
\begin{equation}\label{sage}
\lambda^3-\left(\frac{5}{2}-s\cdot\nabla^2F\cdot\right)\lambda^2+\left(2-\frac{s}{2}\cdot\nabla^2F\right)\lambda-\frac{1}{2}=0.
\end{equation}

According to the stability theory, the absolute stable region is the set of step size $s$ such that all the roots of the characteristic equation (eigenvalues) lie in the unit circle \citep{abs}. Since the left hand of (\ref{sage}) can be factorized to
\[\left(\lambda-\frac{1}{2}\right)\left(\lambda^2-(2-s\cdot\nabla^2F)\lambda+1\right),\] 
this condition equals to $0\leq s\cdot\nabla^2F\leq4$. In another word, the absolutely stable region of SAG is $\left\{s:\;s\cdot\nabla^2F\in[0,4]\right\}$.

We then consider NAG. Applying the same approximations, we get the characteristic equation of NAG
\begin{equation}\label{nage}
\lambda^2-(2-2s\cdot\nabla^2 F)\lambda+(1-s\cdot\nabla^2F)=0.
\end{equation}
It is easy to verify that the absolutely stable region of NAG is $\left\{s:\;s\cdot\nabla^2F\in[0,\frac{4}{3}]\right\}$. In summary, we have the following proposition.
\begin{prop}\label{absp}
	Under abovementioned approximations, the characteristic equation of SAG and NAG are presented in (\ref{sage}) and (\ref{nage}). The absolutely stable region of SAG is $\left\{s:\;s\cdot\nabla^2F\in[0,4]\right\}$, while the absolutely stable region of NAG is $\left\{s:\;s\cdot\nabla^2F\in[0,\frac{4}{3}]\right\}$.
\end{prop}
\begin{rmk}
	If we change the choice of parameters $k, m_1, m_2$ in subsection \ref{nmtd} and so change the form of the new method, the absolutely stable region will remain unchanged.
\end{rmk}
When $s$ lies in the absolutely stable region, stability theory guarantees that errors will not be magnified during iteration. Proposition \ref{absp} shows that the absolutely stable region of SAG is larger than NAG, so our new method is more stable.

\section{Application to matrix completion problem: SFISTA}
Our SAG has a larger absolutely stable region than NAG, so we can choose larger step size when applying SAG and the computation will be faster. In this section we apply our SAG to matrix completion problem. We combine our SAG with proximal operator to present a new algorithm which can be viewed as a modification of the well-known fast iterative shrinkage-thresholding algorithm (FISTA) \citep{fista}. Then we show the better stability and higher computation speed of our new method by comparing with FISTA.

For matrix completion problem there exists a ``true" low rank matrix $\bm{M}$. We are given some entries of $\bm{M}$ and asked to fill missing entries. In our simulated dataset, the size of $\bm{M}$ is $1000\times1000$, and the true rank of $\bm{M}$ is set to be $4$. There have been various algorithms to solve such problem \citep{emm,rag}. Besides, it is proposed that matrix completion can be transformed to the following unconstrained optimization problem \citep{unc}
\begin{equation*}
\min_{\bm{X}}F(\bm{X})=\frac{1}{2}\|\bm{X}_{obs}-\bm{M}_{obs}\|^2+\lambda\|\bm{X}\|_*.
\end{equation*}
$F(\bm{X})$ is composed of a smooth term and a non-smooth term, so gradient-based algorithms cannot be used directly. Proximal gradient algorithms \citep{prox} are widely used in such composite optimization problems, and fast iterative shrinkage-thresholding algorithm (FISTA)  is a successful algorithm. Moreover, FISTA has been extended to matrix completion case \citep{ji}. For convenience, we set $G(\bm{X})=\frac{1}{2}\|\bm{X}_{obs}-\bm{M}_{obs}\|^2$, $H(\bm{X})=\lambda\|\bm{X}\|_*$, and $g(\bm{X})=\nabla G(\bm{X}).$

The idea of FISTA builds on NAG. We also apply accelerated proximal gradient method (APG) \citep{prox} for our numerical experiment, which is composed of NAG and proximal gradient descent. These two algorithms are presented in Appendix \ref{a5}. We find the performances of them are similar in our experiments. 

Our contribution is the third method (Algorithm \ref{sfi}), SAG combined with proximal operator, which we call SFISTA.
\begin{algorithm}
	\caption{SFISTA}
	\label{sfi}
	\begin{algorithmic}
		\STATE \textbf{Input: }step size $s$.
		\STATE \textbf{Initial value: }$\bm{X}_2=\bm{X}_1=\bm{X}_0=\bm{M}_{obs}$.
		\STATE \textbf{$\bm{(k-1)}$th iteration $(k\geq2)$.} Compute
		\begin{align*}
		&\bm{Y}_k=\frac{10k^2+9k+6}{4k^2+8k}\bm{X}_k-\frac{4k^2+3}{2k^2+4k}\bm{X}_{k-1}+\frac{2k-1}{4k+8}\bm{X}_{k-2},\\
		&\bm{Z}_k=\frac{2k-3}{k}\bm{X}_k-\frac{k-3}{k}\bm{X}_{k-1},\\
		&\bm{X}_{k+1}=\mathop{\arg\min}_{\bm{X}}\left\{\frac{1}{2}\cdot\frac{2k+4}{ks}\left\|\bm{X}-\left(\bm{Y}_k-\frac{ks}{2k+4}g(\bm{Z}_k)\right)\right\|^2+\lambda\|\bm{X}\|_*\right\}.
		\end{align*}
	\end{algorithmic}
\end{algorithm}
Notice that the minimizing problems in iterations of above three algorithms can be solved directly by singular value decomposition \citep{svd}. 

\begin{figure}[ht]
	\centering
	\subfigure[original]{
		\includegraphics[width=4.3cm]{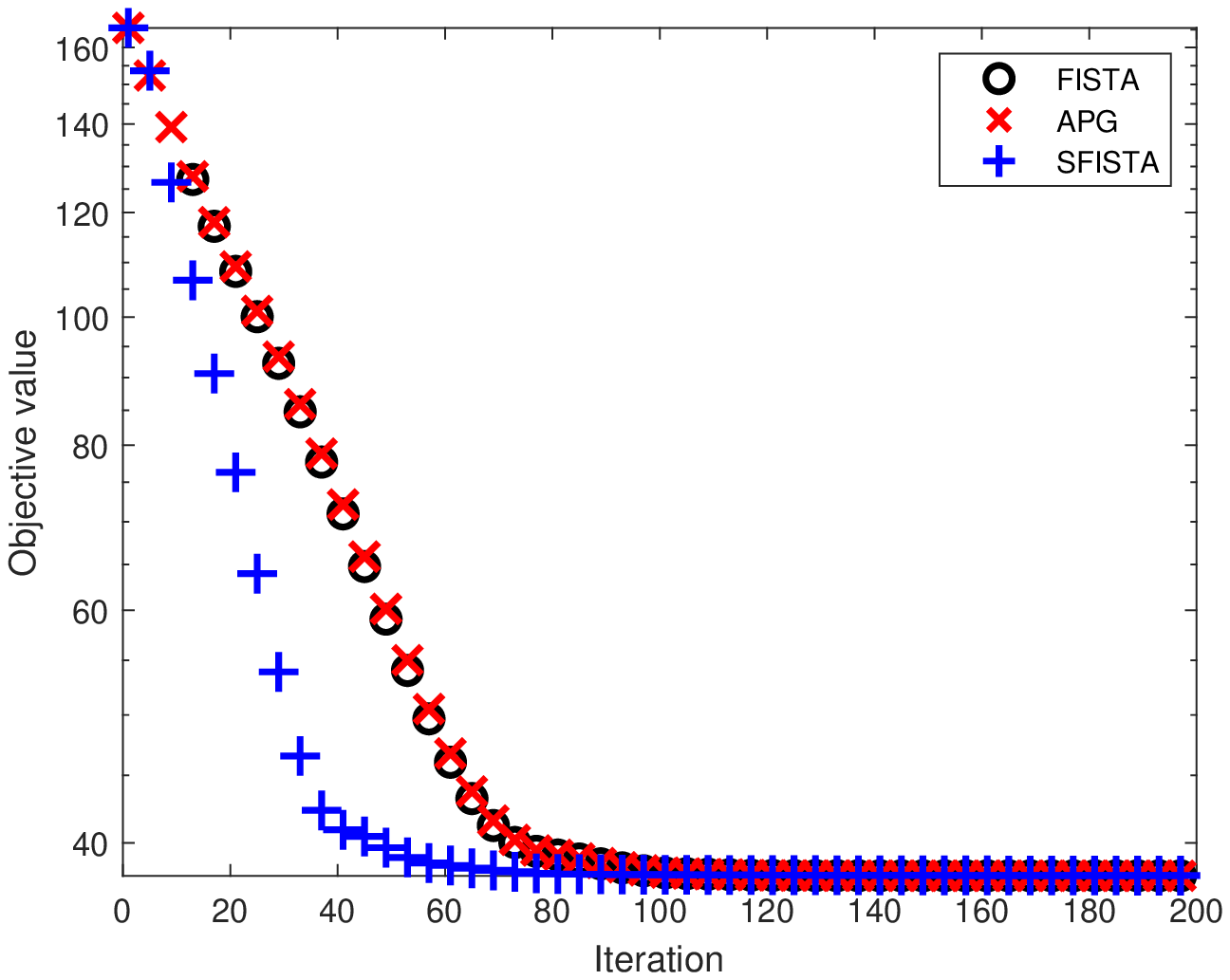}
		\label{f1_o}
	}
	\subfigure[reduced to $80\%$]{
		\includegraphics[width=4.3cm]{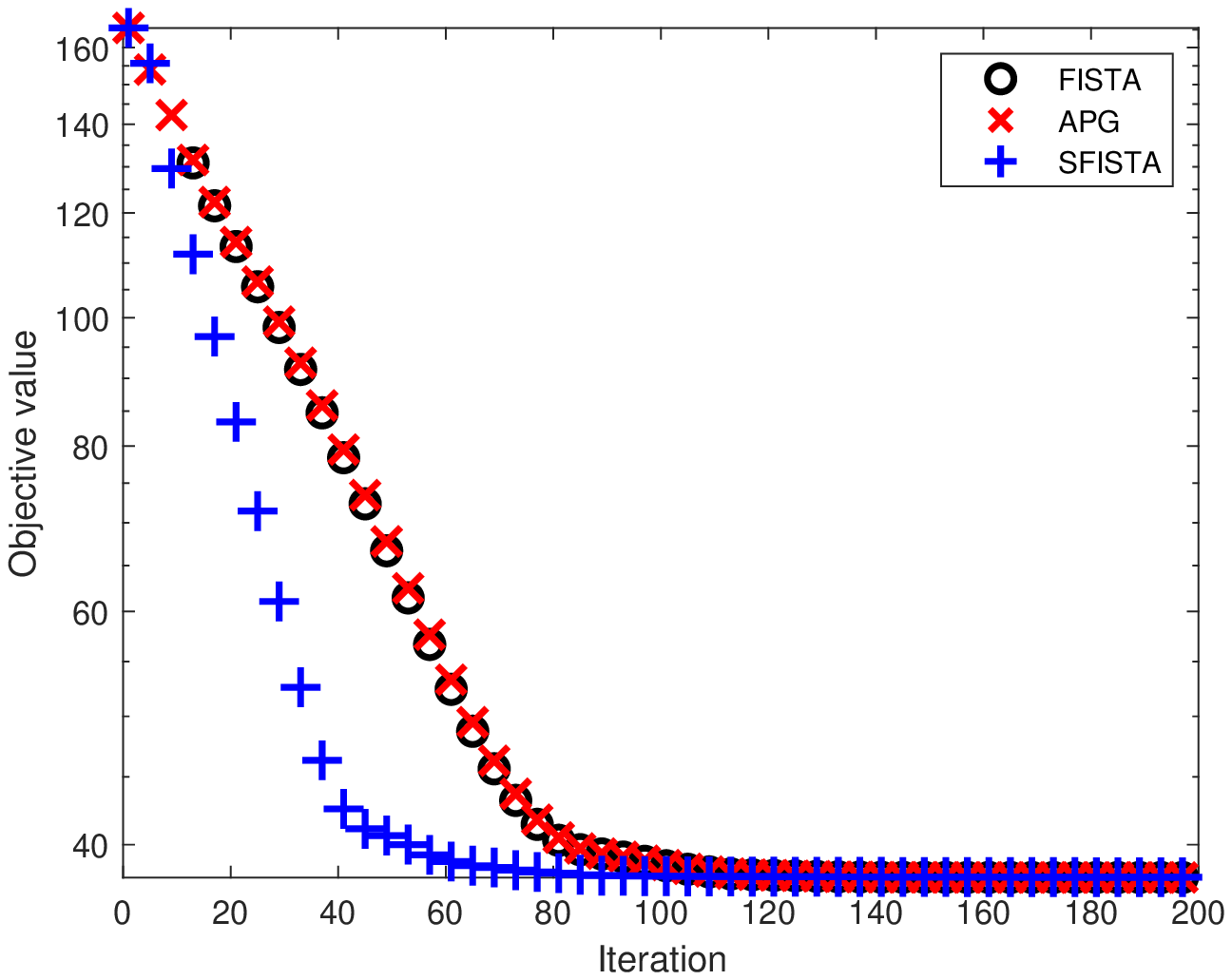}
		\label{f1_8}
	}
	\subfigure[reduced to $50\%$]{
		\includegraphics[width=4.3cm]{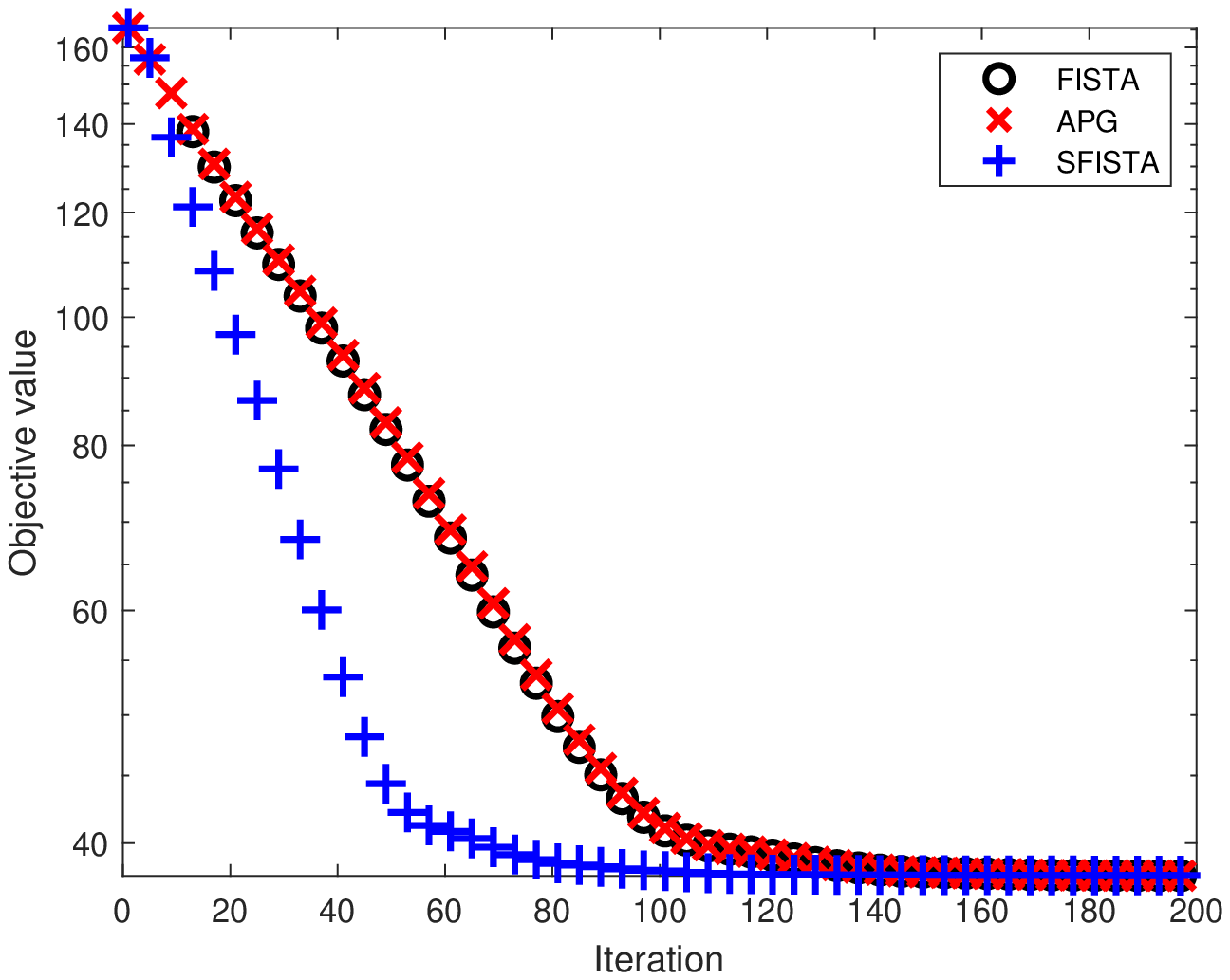}
		\label{f1_5}
	}\\
	\subfigure[reduced to $10\%$]{
		\includegraphics[width=4.3cm]{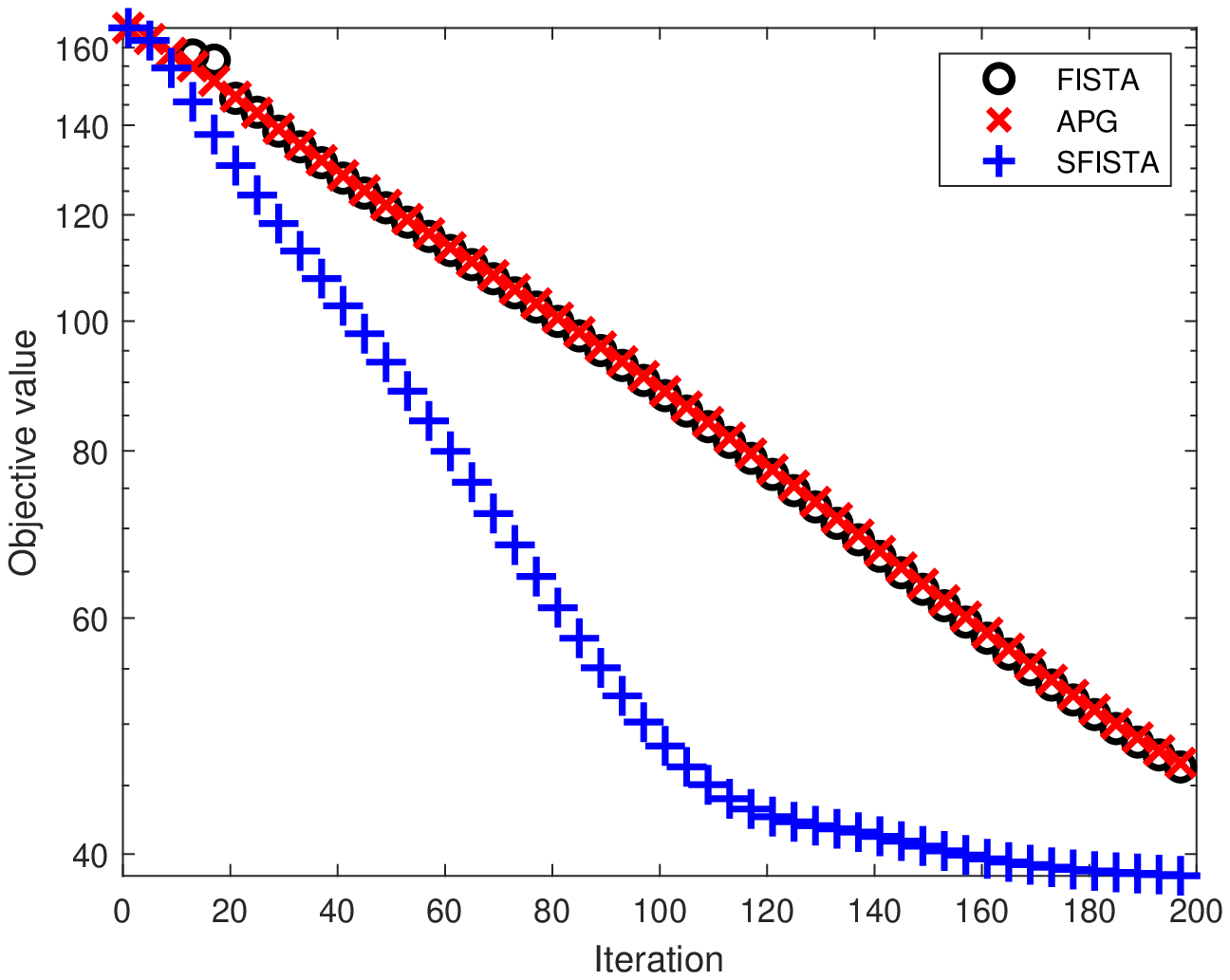}
		\label{f1_1}
	}
	\subfigure[backtracking]{
		\includegraphics[width=4.3cm]{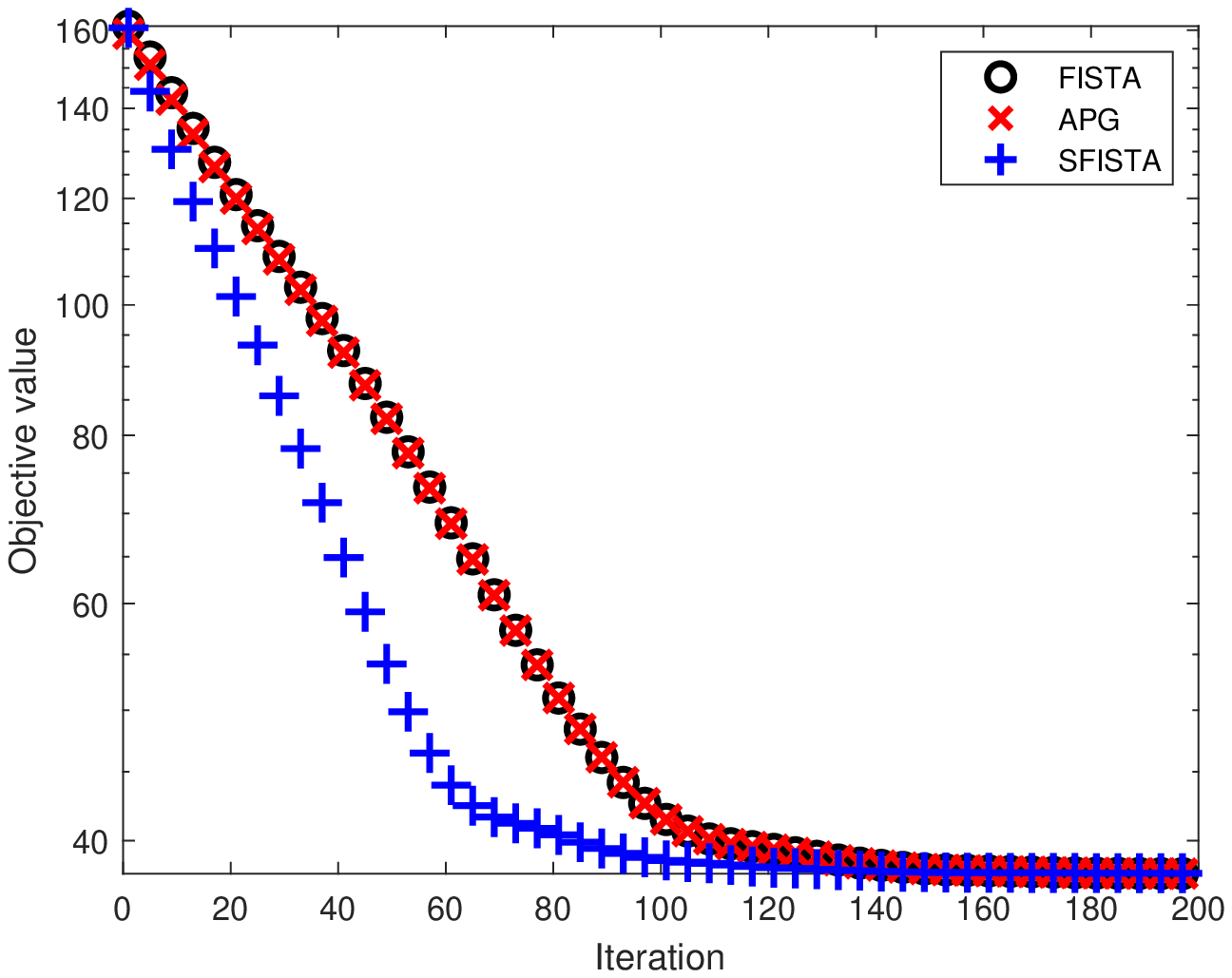}
		\label{bc}
	}
	\caption{Iterations of FISTA, APG and SFISTA for matrix completion objective function. Y-axis represents the objective function $F(x_n)$. In Figure \ref{f1_o}, step size is $1.4$ for FISTA and APG, and $4.5$ for SFISTA. In the next three figures, step sizes are reduced from $1.4$ and $4.5$ in the proportion marked below the figures. In Figure \ref{bc}, step sizes are chosen automatically by backtracking.}
	\label{ffix}
\end{figure}

In the above three algorithms, we use fixed step sizes, and the results are presented in the first four subfigures of Figure \ref{ffix}. During this experiment, our first find is that the range of step size that our SFISTA can converge stably is significantly larger than those of FISTA and APG. Experiment shows that the feasible regions of step size of FISTA and APG are both $(0,1.4]$ (accurate to one decimal place), while that of our SFISTA is $(0,4.5]$. Secondly, we find empirically that for all methods, convergence is faster when step size is larger (within the feasible region). Therefore, we choose the largest step sizes in the feasible region for all methods to compare their highest computational speed. We also compare performances of the three methods with step sizes reduced from the largest in equal proportion. We find that when step sizes are chosen to be the largest or reduced from the largest in equal proportion ($80\%$, $50\%$, $10\%$), convergence of SFISTA is faster than the other two methods.

Furthermore, we combine the three methods with backtracking \citep{fista} to choose step sizes automatically. We present SFISTA with backtracking below (Algorithm \ref{mfls}), and the other two algorithms are similar (presented in Appendix \ref{a5}). 
\begin{algorithm}[ht]
	\caption{SFISTA with backtracking}
	\label{mfls}
	\begin{algorithmic}
		\STATE \textbf{Input: }some $\beta<1$.
		\STATE \textbf{Initial value: }$\bm{X}_2=\bm{X}_1=\bm{X}_0=\bm{M}_{obs}$, step size $s_2$.
		\STATE $\bm{(k-1)}$\textbf{th iteration $(k\geq2)$.}
		\begin{align*}
		&\bm{Y}_k=\frac{10k^2+9k+6}{4k^2+8k}\bm{X}_k-\frac{4k^2+3}{2k^2+4k}\bm{X}_{k-1}+\frac{2k-1}{4k+8}\bm{X}_{k-2},\\
		&\bm{Z}_k=\frac{2k-3}{k}\bm{X}_k-\frac{k-3}{k}\bm{X}_{k-1}.
		\end{align*}
		\STATE \hspace*{1em}Find the smallest nonnegative integer $i_{k+1}$ such that with $s=\beta^{i_{k+1}}s_k$
		\[F(\widetilde{\bm{X}})<F(\bm{Y}_k)+\left\langle\widetilde{\bm{X}}-\bm{Y}_k,g(\bm{Z}_k)\right\rangle+\frac{1}{2}\cdot\frac{2k+4}{ks}\|\widetilde{\bm{X}}-\bm{Y}_k\|^2,\]
		\hspace*{1em}where
		\[\widetilde{\bm{X}}=\mathop{\arg\min}_{\bm{X}}\left\{\frac{1}{2}\cdot\frac{2k+4}{ks}\left\|\bm{X}-\left(\bm{Y}_k-\frac{ks}{2k+4}g(\bm{Z}_k)\right)\right\|^2+\lambda\|\bm{X}\|_*\right\}.\]
		\STATE \hspace*{1em}Set $s_{k+1}=\beta^{i_{k+1}}s_k$ and compute
		\begin{align*}
		\bm{X}_{k+1}=\widetilde{\bm{X}}.
		\end{align*}
	\end{algorithmic}
\end{algorithm}
Performances of the three algorithms with backtracking are compared on abovementioned dataset, and the result is shown in Figure \ref{bc}. Convergence of SFISTA is faster than the other two methods. Moreover, we find that the step size of SAG is reduced $9$ times during the $200$ iterations, while step sizes of FISTA and APG are both reduced $13$ times. This find also confirms the better stability of SFISTA.

In short, the better stability of SFISTA enables it to work well with larger step size than FISTA and APG, which also leads to faster computation.

\section{Application to handwriting digit recognition}
In this section, we evaluate SAG and NAG with handwriting digit recognition task. \emph{MNIST} \citep{mni} is a large dataset of handwriting digits and is widely used in training and testing image recognition system. One of the most popular models to analyze visual imagery is convolutional neural network (CNN) \citep{lecun}. Here we train a two-layer CNN by SAG and NAG on \emph{MNIST} dataset and evaluate their performances.

\begin{figure}[ht]
	\centering
	\subfigure[$s=0.02$]{
		\includegraphics[width=4.3cm]{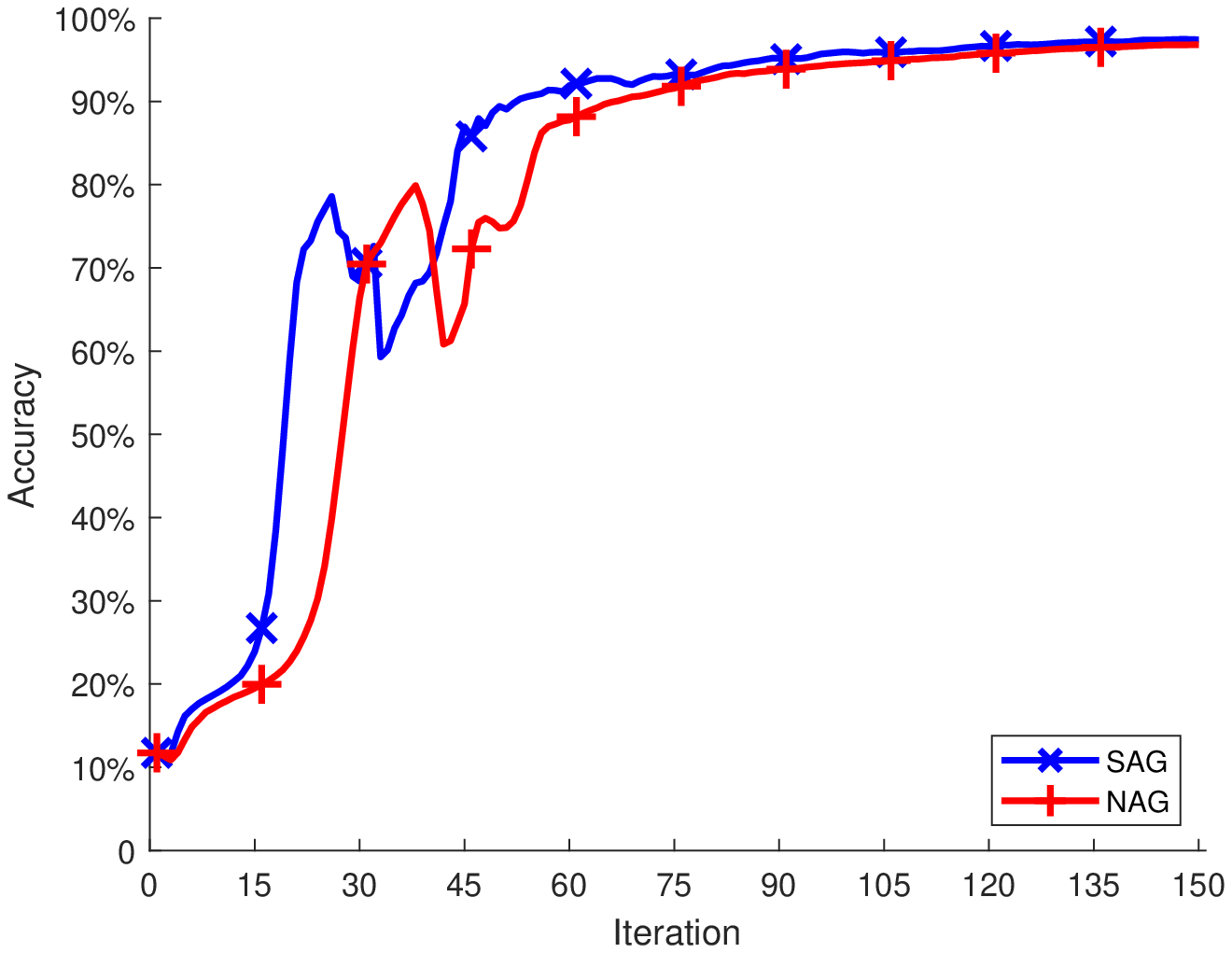}
		\label{f2_1}
	}
	\subfigure[$s=0.08$]{
		\includegraphics[width=4.3cm]{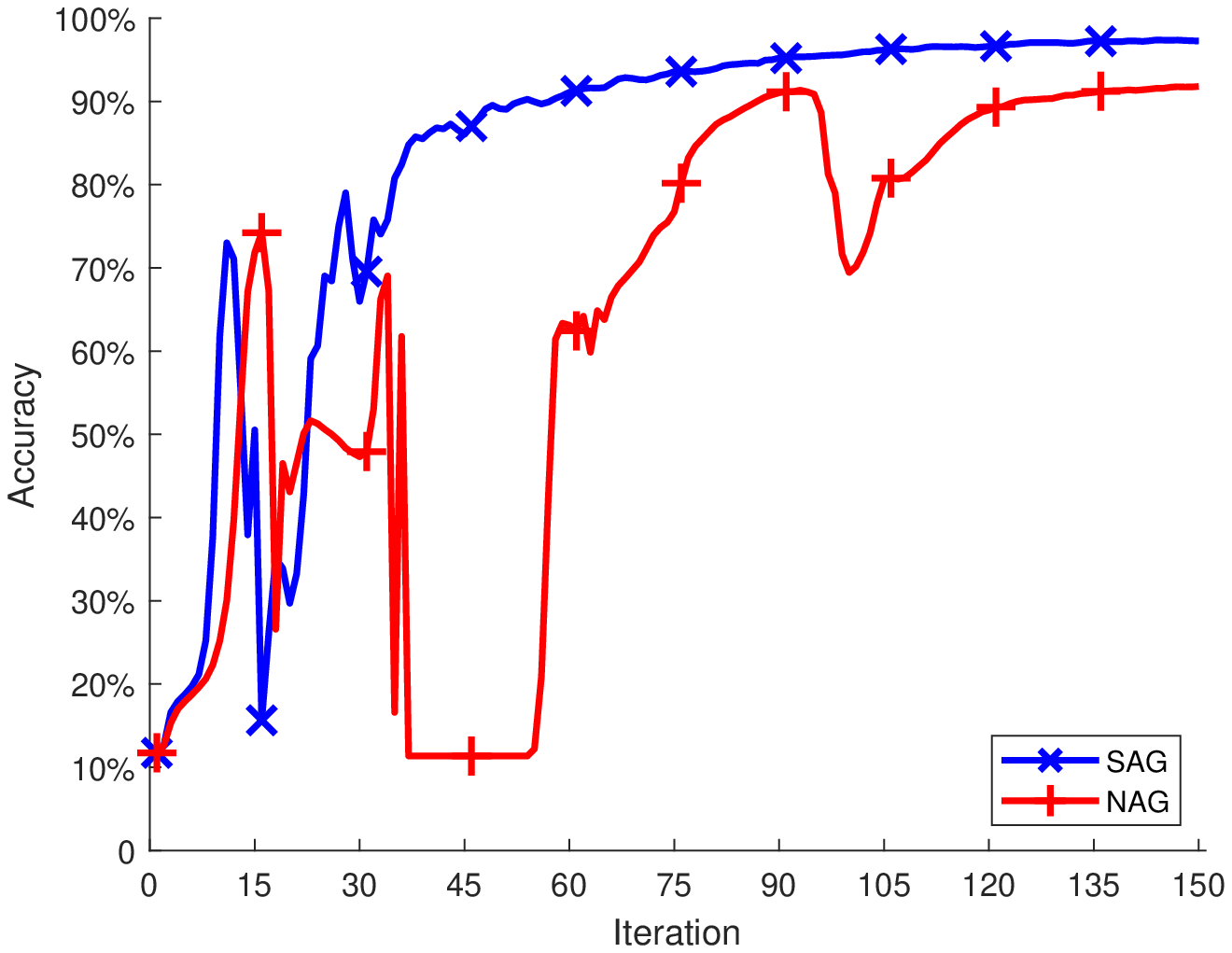}
		\label{f2_2}
	}
	\subfigure[$s=0.14$]{
		\includegraphics[width=4.3cm]{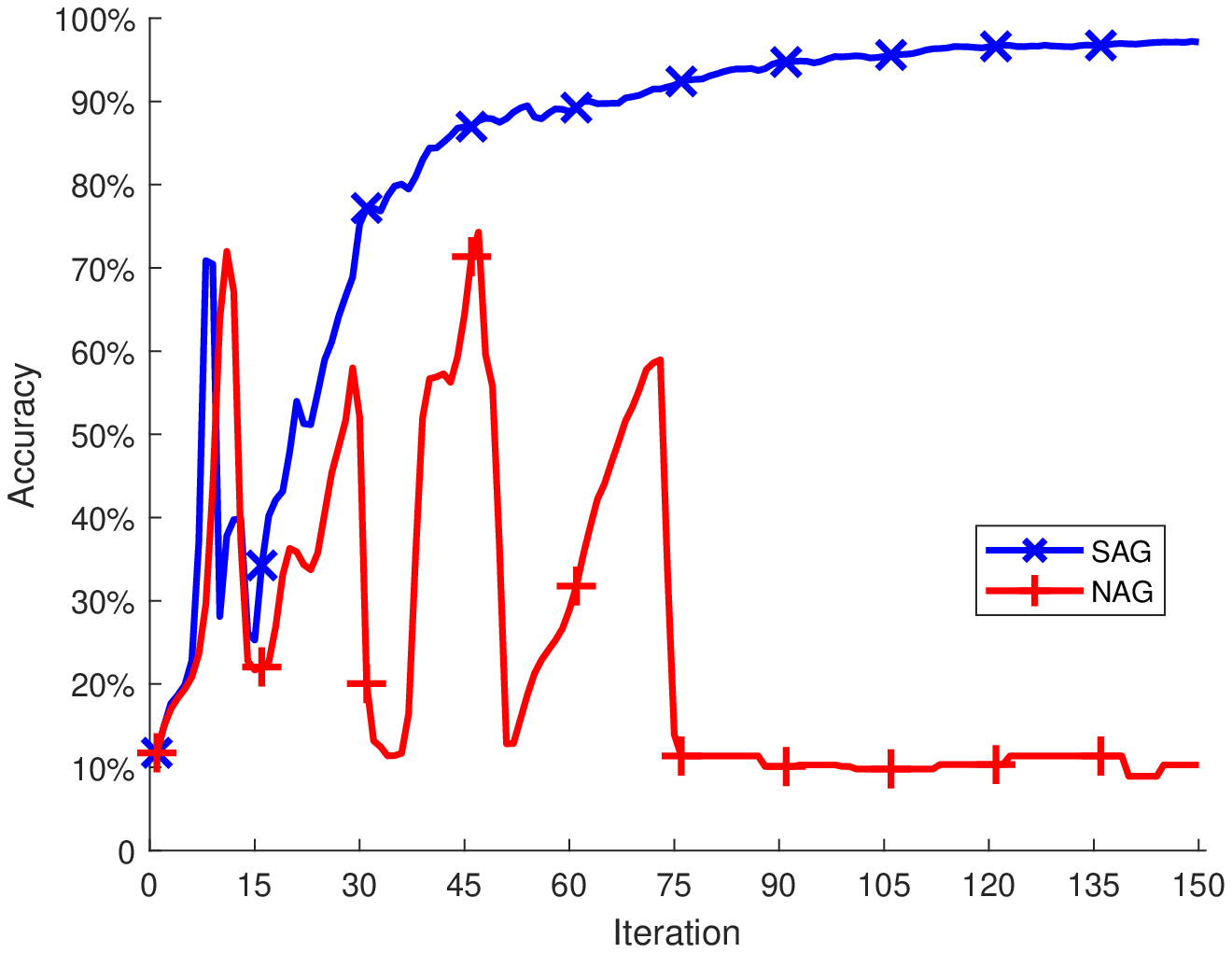}
		\label{f2_3}
	}
	\caption{Iterations of SAG and NAG for test accuracy. Y-axis represents the test accuracies for both methods. Step size is $0.02$ in Figure \ref{f2_1}, $0.08$ in Figure \ref{f2_2} and $0.14$ in Figure \ref{f2_3}.}
	\label{fmni}
\end{figure}

\emph{MNIST} dataset \citep{mni} contains $60000$ training images and $10000$ test images. In our experiment, we use the first $20000$ images to train a two-layer convolutional neural network and use the whole $10000$ images to test the accuracy of our network. The test accuracies during iteration are provided in Figure \ref{fmni}. When the step size is relatively small ($0.02$ in the experiment), accuracies of both methods increase to more than $97\%$ at similar speeds. When the step size is enlarged to $0.08$, the accuracy curve of NAG is relatively unstable compared to our SAG. When the step size is large ($0.14$ in the experiment), accuracy of SAG can still increase to $97\%$, while NAG oscillates violently and fails to train a useful network. We find that our SAG has less oscillation than NAG, which also leads to larger feasible region of step size and faster convergence. This consists with our theoretical result and confirms the better stability of SAG. More detailed results of this experiment are presented in Appendix \ref{a6}.

To compare SAG and NAG more quantitatively, we take another experiment to evaluate the local performance of the two methods. Note that we use quadratic approximation in stability analysis. Therefore, due to the non-convexity of neural network objective function, the better performance of SAG is more significant when parameters of the network are not far from the optimal ones. More precisely, our training procedure consists of two steps:
\begin{enumerate}
	\item Train the neural network with an approximate method until the test accuracy is not less than $90\%$,
	\item Switch to optimizer SAG/NAG with step size $s$.
\end{enumerate}
Then we compare the performance of SAG and NAG with different step size in step 2.

In the experiment, we randomly choose $5$ different initial values and use NAG with sufficiently small step size $(0.02)$ for step 1. For step 2, we apply SAG and NAG with $30$ different step sizes (from $0.01$ to $0.3$ with gap $0.01$), and measure the performance by the number of steps needed to gain a loss value less than $0.1$. In Table \ref{t_step}, we record the best choices of step size for different initial values and algorithms. It is obvious the best step size of SAG is about two times as large as that of NAG. We also present the number of steps needed with the best step size for different initial values and algorithms in table \ref{t_num}. The results show that our SAG is about $1.6$ times faster than NAG in step 2.
\begin{table}[ht]
	\centering
	\begin{tabular}{*{6}{c}}
		\toprule
		Algorithm & Init.1 & Init.2 & Init.3 & Init. 4 & Init. 5\\
		\midrule
		SAG & \textbf{0.26} & \textbf{0.19} & \textbf{0.07} & \textbf{0.24} & \textbf{0.15}\\
		NAG & 0.13 & 0.1 & 0.04 & 0.11 & 0.08\\
		\bottomrule
	\end{tabular}
	\caption{Best choice of step size for SAG and NAG.}
	\label{t_step}
\end{table}
\begin{table}[ht]
	\centering
	\begin{tabular}{*{6}{c}}
		\toprule
		Algorithm & Init.1 & Init.2 & Init.3 & Init. 4 & Init. 5\\
		\midrule
		SAG & \textbf{28} & \textbf{32} & \textbf{41} & \textbf{27} & \textbf{31}\\
		NAG & 48 & 53 & 56 & 51 & 46\\
		\bottomrule
	\end{tabular}
	\caption{Minimum number of steps needed to achieve $0.1$ loss for SAG and NAG.}
	\label{t_num}
\end{table}

\section{Discussion}\label{s6}
In this paper we give the precise order of NAG as a numerical approximation of the limit differential equation. Inspired by this perspective, we present a new method which we call \emph{stabilized accelerated gradient} (SAG). 

We prove that our SAG has a larger absolutely stable region, so it works more stable with large step size. We evaluate the performance of SAG in two experiments: matrix completion and handwriting digit recognition. In both experiments, we find that our SAG has a larger region of feasible step size than NAG, which confirms the better stability of SAG. Moreover, better stability of SAG leads to higher computational speed.

In summary, we provide a new way to improve the stability of Nesterov's accelerated gradient method, and our approach indicates how properties in numerical analysis can be used to refine an optimization method. We believe this approach is worth further study and more higher order methods may be worth considering. Finally, more theoretical analysis of SAG is still a future topic. From previous paper and this work we know that the limit ODE (\ref{ode}) converges to the minimum at a speed $\mathcal{O}(t^{-2})$ \citep{su} (coincides with $\mathcal{O}(n^{-2})$ of NAG) and SAG is a better approximation of the ODE than Nesterov's method. These two facts imply the $\mathcal{O}(n^{-2})$ convergence of SAG, which is also verified by the better performance of SAG in practice. However, a direct convergence analysis is to be involved. This analysis might be finished by technical construction of estimate sequence. Furthermore, it would be interesting to apply SAG to deep neural network, especially to a landscape with large gradient.

\bibliographystyle{plainnat}
\bibliography{nips_ref}

\appendix

\section{Proof of Theorem \ref{trun}}\label{a1}
Theorem \ref{trun}. Assume $\nabla F$ satisfies $L$-Lipschitz condition, and solution $x(t)$ of the derived differential equation (\ref{ode}) has a continuous third derivative. For fixed time $t$, the truncation error (\ref{rn}) satisfies
\begin{equation}\label{ro1}
L[x(t);h]=\mathcal{O}(h^3).
\end{equation}
\begin{proof}
	We substitute the equality
	\[x(t-h)=x(t)+\mathcal{O}(h)\]
	to the last term of
	\begin{align}\label{rrn}
	\begin{split}
	L[x(t);h]=&x(t+h)-\frac{2t-3h}{t}x(t)+\frac{t-3h}{t}x(t-h)+\\&h^2\nabla F\left(x(t)+\frac{t-3h}{t}\left(x(t)-x(t-h)\right)\right)
	\end{split}
	\end{align}
	to get
	\begin{align*}
	\nabla F\left(x(t)+\frac{t-3h}{t}(x(t)-x(t-h))\right)&=\nabla F\left(x(t)+\frac{t-3h}{t}\cdot \mathcal{O}(h)\right).
	\end{align*}
	Since $\nabla F$ satisfies $L$-Lipschitz condition, we know
	\begin{align*}
	\nabla F\left(x(t)+\frac{t-3h}{t}(x(t)-x(t-h))\right)=&\nabla F(x(t))+\mathcal{O}(h)\\
	=&-x^{(2)}(t)-\frac{3}{t}x^{(1)}({t})+\mathcal{O}(h).
	\end{align*}
	For the second equality, we substitute the differential equation (\ref{ode}). Then we expend the first and third terms of $L[x(t);h]$ to third order
	\begin{align*}
	x(t+h)=x(t)+hx^{(1)}(t)+\frac{h^2}{2}x^{(2)}(t)+\mathcal{O}(h^3),\\
	x(t-h)=x(t)-hx^{(1)}(t)+\frac{h^2}{2}x^{(2)}(t)+\mathcal{O}(h^3).
	\end{align*}
	Finally, we substitute these three equations to the truncation error (\ref{rrn}) to conclude
	\[L[x(t);h]=\mathcal{O}(h^3).\qedhere\]
\end{proof}
\begin{rmk}
	(\ref{ro1}) can be written as 
	\begin{equation*}
	|L[x(t);h]|\leq M_0h^3,
	\end{equation*}
	where $M_0$ depends on $\sup_{s\leq t}|x^{(1)}(s)|$ and $\sup_{s\leq t}|x^{(3)}(s)|$. 
\end{rmk}
\begin{rmk}\label{r3}
	Theorem \ref{trun} deals with the problem for fixed time $t$. To finish the proof of the approximation theorem, we have to consider the situation that $t_n=nh$, where $n\geq1$ is fixed.
	
	We set a fixed time $t_0$ and assume that $t_n=nh<t_0$. Since $x(t)$ has a continuous third derivative, $x(t)$ and its first to third derivative are bounded in $[0,t_0]$. We replace time $t$ in the above proof by $t_n$ and expend the terms of (\ref{rrn}). Now the term
	\[-\frac{3h^3}{2t_n}x^{(2)}(t_n)\]
	obtained from the expansion of $x(t_n-h)$ cannot be viewed as $\mathcal{O}(h^3)$, but there exists $M_2>0$ such that 
	\[\left|-\frac{3h^3}{2t_n}x^{(2)}(t_n)\right|\leq M_2\frac{h^2}{n}.\]
	As a consequence, we have 
	\begin{equation}\label{lmi}
	|L[x(t_n);h]|\leq M_1h^3+M_2\frac{h^2}{n},
	\end{equation}
	where $M_1$ and $M_2$ rely on $t_0$.
\end{rmk}
\section{Proof of Lemma \ref{mai}}\label{a2}
	Lemma \ref{mai}. We define matrices $\bm{C}_n$ and $\bm{D}_{n,l}$ as
\begin{align*}
&\bm{C}_n=\begin{pmatrix}\frac{2n-1}{n+1}&&-\frac{n-2}{n+1}\\1&&0\end{pmatrix},\\
&\bm{D}_{n,l}=\bm{C}_n\bm{C}_{n-1}\cdots\bm{C}_{n-l+1},
\end{align*}
where $n\geq0$ and $0<l\leq n+1$. In addition, we set $\bm{D}_{n,0}=\bm{I}_2$. Then there exist positive constants $M,\;M_3$ such that for all $n$, the following two inequalities hold, where the matrix norm is 2-norm.
\begin{align}
\sup_{0\leq l\leq n+1}\|\bm{D}_{n,l}\|&\leq Mn,\label{m}\\
\bm{D}_{n,n+1}&\leq M_3.\label{mt}
\end{align}
\begin{proof}
	From
	\[\bm{C}_2=\begin{pmatrix}1&&0\\1&&0\end{pmatrix},\]
	we have 
	\[\bm{D}_{n,n-1}=\begin{pmatrix}1&&0\\1&&0\end{pmatrix},\quad\bm{D}_{n,n}=\begin{pmatrix}\frac{1}{2}&&\frac{1}{2}\\\frac{1}{2}&&\frac{1}{2}\end{pmatrix},\quad\bm{D}_{n,n+1}=\begin{pmatrix}0&&1\\0&&1\end{pmatrix},\]
	for $n\geq2$. So it is obvious that there exists $M_3$ to make (\ref{mt}) true and $M_4>0$ such that for all $n<2$ or $n\geq2,\;l>n-2$ or $l=0$, \begin{equation}\label{m0}
	\|\bm{D}_{n,l}\|\leq M_4n.
	\end{equation}
	
	Then we consider the situation of $n\geq2,\;0<l\leq n-2$. Notice that 
	\[\bm{C}_n=\bm{P}\begin{pmatrix}1&&1\\0&&\frac{n-2}{n+1}\end{pmatrix}\bm{P}^{-1},\]
	where
	\[\bm{P}=\begin{pmatrix}1&&1\\1&&0\end{pmatrix}.\]
	Assume we have already got
	\[\bm{D}_{n,l}=\bm{P}\begin{pmatrix}1&&a_{n,l}\\0&&b_{n,l}\end{pmatrix}\bm{P}^{-1}\]
	satisfying
	\[0<a_{n,l}\leq l,\qquad0<b_{n,l}\leq1.\]
	Then since
	\[\bm{D}_{n,l+1}=\bm{D}_{n,l}\bm{C}_{n-l}\qquad and\qquad0\leq\frac{n-l-2}{n-l+1}<1,\]
	$\bm{D}_{n,l+1}$ can be written as
	\[\bm{D}_{n,l+1}=\bm{P}\begin{pmatrix}1&&a_{n,l+1}\\0&&b_{n,l+1}\end{pmatrix}\bm{P}^{-1},\]
	where
	\[0<a_{n,l+1}\leq l+1,\qquad0<b_{n,l}\leq1.\]
	Then for fixed $n$, we induce from $l=1$ to get
	\begin{equation*}
	\bm{D}_{n,l}=\bm{P}\widetilde{\bm{D}}_{n,l}\bm{P}^{-1}\triangleq\bm{P}\begin{pmatrix}1&&a_{n,l}\\0&&b_{n,l}\end{pmatrix}\bm{P}^{-1},
	\end{equation*}
	where
	\begin{equation}\label{ctrl}
	0<a_{n,l}\leq l\leq n,\qquad0<b_{n,l}\leq1,
	\end{equation}
	for all $n\geq2,\;0<l\leq n-2$. Then we can estimate $\|\bm{D}_{n,l}\|$. Calculation shows that
	\[\widetilde{\bm{D}}_{n,l}\widetilde{\bm{D}}_{n,l}^\mathsf{T}=\begin{pmatrix}1+a_{n,l}^2&&a_{n,l}b_{n,l}\\a_{n,l}b_{n,l}&&a_{n,l}^2\end{pmatrix}.\]
	The eigenvalues of this matrix are
	\[\lambda_{1,2}=\frac{1+a_{n,l}^2+b_{n,l}^2\pm\sqrt{(1+a_{n,l}^2+b_{n,l}^2)^2-4b^4}}{2}.\]
	Combining this representation with (\ref{ctrl}), we get the estimation
	\begin{align*}
	\|\widetilde{\bm{D}}_{n,l}\|=\sqrt{\max\{|\lambda_1|,|\lambda_2|\}}\leq\sqrt{1+a_{n,l}^2+b_{n,l}^2}\leq n+2.
	\end{align*}
	So there exists $M_5>0$, such that for all $n\geq2,\;0<l\leq n-2$, inequality
	\begin{equation}\label{m1}
	\|\bm{D}_{n,l}\|\leq M_5n
	\end{equation}
	holds. Combining (\ref{m0}) with (\ref{m1}), we finish the proof of (\ref{m}). 
\end{proof}
\section{Proof of Theorem \ref{conv}}\label{a3}
Theorem \ref{conv}. Under conditions in Theorem \ref{trun}, for fixed time $t$ and $n=t/h$, $x_n$ converges to $x(t)$ at a rate of $\mathcal{O}(h\ln\frac{1}{h})$ if $x_0=x(0)$ and $x_1=x(h)$.
\begin{proof}
	In this proof, we first calculate the error caused by a single iteration, which can be divided into an accumulation term and a truncation term. Then we use the error estimation given by Theorem \ref{trun} and apply discrete Gronwall inequality to prove the convergence.
	
	The difference of our situation from classic situation in numerical analysis is that the iteration scheme changes with $n$, so the classic technique can not be used to bound the norm of the transition matrix and the truncation error. Our new approach is presented in Remark \ref{r3} of Theorem \ref{trun} and Lemma \ref{mai}, and is the most important technical innovation in our proof.
	
	Recall the recurrence relation
	\[x_{n+1}=x_n+\frac{n-3}{n}(x_n-x_{n-1})-h^2\nabla F\left(x_n+\frac{n-3}{n}(x_n-x_{n-1})\right)\]
	and the definition of truncation error
	\[x(t_{n+1})=x(t_n)+\frac{n-3}{n}(x(t_n)-x(t_{n-1}))-h^2\nabla F\left(x(t_n)+\frac{n-3}{n}(x(t_n)-x(t_{n-1}))\right)+L[x(t_n);h],\]
	where $t_n=nh$.
	
	We subtract the above two equations, and introduce overall error
	\[e_n=x(t_n)-x_n\]
	to get
	\[e_{n+1}=\frac{2n-3}{n}e_n-\frac{n-3}{n}e_{n-1}-h^2b_{n-1}+L[x(t_n);h],\]
	which can also be written as
	\begin{equation}\label{err}
	e_{n+2}-\frac{2n-1}{n+1}e_{n+1}+\frac{n-2}{n+1}e_{n}=-h^2b_{n}+L[x(t_{n+1});h],
	\end{equation}
	where 
	\begin{equation}\label{bn}
	b_n=\nabla F\left(\frac{2n-1}{n+1}x_{n+1}-\frac{n-2}{n+1}x_n\right)-\nabla F\left(\frac{2n-1}{n+1}x(t_{n+1})-\frac{n-2}{n+1}x(t_n)\right).
	\end{equation}
	Then we rewrite (\ref{err}) into a form that is convenient for recurrence. We set
	\[\bm{E}_n=\begin{pmatrix}e_{n+1}\\e_n\end{pmatrix},\quad\bm{C}_n=\begin{pmatrix}\frac{2n-1}{n+1}&&-\frac{n-2}{n+1}\\1&&0\end{pmatrix},\quad\bm{B}_n=\begin{pmatrix}-h^2b_n^*\\0\end{pmatrix},\]
	where
	\[b_n^*=-\frac{e_{n+2}-\frac{2n-1}{n+1}e_{n+1}+\frac{n-2}{n+1}e_{n}}{h^2}=b_n-\frac{L[x(t_{n+1});h]}{h^2}.\]
	Then (\ref{err}) can be written as
	\begin{equation*}
	\bm{E}_{n+1}=\bm{C}_n\bm{E}_n+\bm{B}_n.
	\end{equation*}
	By recursive method, we have
	\begin{equation*}
	\bm{E}_n=\bm{C}_{n-1}\cdots\bm{C}_0\bm{E}_0+\sum_{l=1}^{n}\bm{C}_{n-1}\cdots\bm{C}_{n-l+1}\bm{B}_{n-l}.
	\end{equation*}
	With the notations introduced in Lemma \ref{mai}, this equality can be written as
	\begin{equation}\label{recur}
	\bm{E}_n=\bm{D}_{n-1,n}\bm{E}_0+\sum_{l=1}^n\bm{D}_{n-1,l-1}\bm{B}_{n-l}.
	\end{equation}
	Now we need to estimate $\|\bm{B}_n\|$. From (\ref{bn}) and the $L$-Lipschitz proporty of $\nabla F$, we have
	\[|b_n|\leq L\left(\frac{2n-1}{n+1}|e_{n+1}|+\frac{n-2}{n+1}|e_n|\right)\leq L\left(2|e_{n+1}|+|e_n|\right)\leq3L\|\bm{E}_n\|\]
	and
	\begin{equation}\label{ectrl}
	\|\bm{B}_n\|\leq 3h^2L\|\bm{E}_n\|+L[x(t_{n+1});h].
	\end{equation}
	Taking norm on both sides of (\ref{recur}) and substituting (\ref{ectrl}) and conclusion of Lemma \ref{mai} yields
	\begin{align}\label{sum0}
	\begin{split}
	\|\bm{E}_n\|&\leq M_3\|\bm{E}_0\|+M(n-1)\sum_{l=0}^{n-1}\left(3h^2L\|\bm{E}_l\|+L[x(t_{l+1});h]\right)\\
	&\leq M_3\|\bm{E}_0\|+3Mnh^2L\sum_{l=0}^{n-1}\|\bm{E}_l\|+Mn\sum_{l=0}^{n-1}L[x(t_{l+1});h].
	\end{split}
	\end{align}
	Now we deal with truncation errors. Recall (\ref{lmi}) in Remark \ref{r3} of Theorem \ref{trun}
	\[|L[x(t_l);h]|\leq M_1h^3+M_2\frac{h^2}{l}.\]
	Take sum to obtain
	\begin{equation}\label{sum1}
	\sum_{l=0}^{n-1}|L[x(t_{l+1});h]|\leq nM_1h^3+M_2h^2\sum_{l=0}^{n-1}\frac{1}{l+1}.
	\end{equation}
	In addition, we have the classic inequality
	\[\sum_{i=1}^{n}\frac{1}{i}\leq\ln n+M_e,\]
	where $M_e$ refers to a positive constant. We substitute it to (\ref{sum1}) to get
	\[\sum_{l=0}^{n-1}|L[x(t_{l+1});h]|\leq nM_1h^3+M_2h^2(\ln n+M_e).\]
	Then we substitute this inequality to (\ref{sum0}) to get a control of $\|\bm{E}_n\|$
	\[\|\bm{E}_n\|\leq M_3\|\bm{E}_0\|+3Mnh^2L\sum_{l=0}^{n-1}\|\bm{E}_l\|+MM_1n^2h^3+MM_2M_enh^2+MM_2nh^2\ln n.\]
	Using discrete Gronwall inequality, we have
	\begin{align*}
	\begin{split}
	\|\bm{E}_n\|\leq e^{3Mn^2h^2L}\left(M_3\|\bm{E}_0\|+MM_1n^2h^3+MM_2M_enh^2+MM_2nh^2\ln n+3Mnh^2L\|\bm{E}_0\|\right).
	\end{split}
	\end{align*}
	Then for fixed $t$, from the relation $n=\frac{t}{h}$ we get
	\[\|\bm{E}_{t/h}\|\leq e^{3Mt^2L}\left((M_3+3MthL)\|\bm{E}_0\|+(MM_1t^2+MM_2M_et)h+MM_2th\ln\frac{t}{h}\right).\]
	Notice that \[\lim_{h\to0}h\ln\frac{t}{h}=0.\]
	So if $\bm{E}_0=\bm{0}$, then the vector form of overall error $\bm{E}_{t/h}$ satisfies
	\[\lim_{h\to0}\|\bm{E}_{t/h}\|=0.\]
	More precisely, we have
	\[\|\bm{E}_{t/h}\|=\mathcal{O}\left(h\ln\frac{1}{h}\right),\]
	for $h\to0$.
\end{proof}
\section{Proof of Theorem \ref{t_rnorder}}\label{a4}
Theorem \ref{t_rnorder}. If $\nabla F$ has continuous second order derivative, the first and second derivative are bounded, and $x(t)$ has continuous fourth derivative, then for fixed $t$, truncation error of SAG (\ref{ite2e}) satisfies 
\begin{equation*}
L[x(t);h]=\mathcal{O}(h^4).
\end{equation*}

\begin{proof}
	Recall the proof of Theorem \ref{trun}. Now we expand $x(t-h)$ to first order
	\[x(t-h)=x(t)+hx^{(1)}(t)+\mathcal{O}(h^2).\]
	Then we have
	\begin{align}\label{nab1}
	\begin{split}
	&\nabla F\left(x(t)+\frac{t-3h}{t}(x(t)-x(t-h))\right)\\
	=&\nabla F\left(x(t)+\left(1-\frac{3h}{t}\right)(hx^{(1)}(t)+\mathcal{O}(h^2))\right)\\
	=&\nabla F\left(x(t)+hx^{(1)}(t)+\mathcal{O}(h^2)\right)\\
	\overset{(1)}{=}&\nabla F\left(x(t)+hx^{(1)}(t)\right)+\mathcal{O}(h^2)\\
	\overset{(2)}{=}&\nabla F(x(t))+\nabla^2F(x(t))\left(hx^{(1)}(t)\right)+\mathcal{O}(h^2)\\
	=&\nabla F(x(t))+h\cdot\nabla^2F(x(t))\cdot x^{(1)}(t)+\mathcal{O}(h^2).
	\end{split}
	\end{align}
	In the $(1)$ equality, we use the $L$-Lipschitz proporty of $\nabla F$. In the $(2)$ equality, we expand $\nabla F$ to second order. To do this, we need $\nabla F$ has continuous second derivative and the second derivative is bounded. 
	
	Then we take derivative on both sides of differential equation
	\[x^{(2)}(t)+\frac{3}{t}x^{(1)}(t)+\nabla F(x(t))=0\]
	to get 
	\begin{equation}\label{nab2}
	\nabla^2F(x(t))\cdot x^{(1)}(t)=\frac{\df}{\df t}(\nabla F(x(t)))=-x^{(3)}(t)-\frac{3}{t}x^{(2)}(t)+\frac{3}{t^2}x^{(1)}(t).
	\end{equation}
	The differential equaltion can also be written as
	\begin{equation}\label{nab3}
	\nabla F(x(t))=-x^{(2)}(t)-\frac{3}{t}x^{(1)}(t).
	\end{equation}
	We substitute (\ref{nab2}) and (\ref{nab3}) to (\ref{nab1}) to get
	\begin{align}\label{fexp}
	\begin{split}
	\nabla F\left(x(t)+\frac{t-3h}{t}(x(t)-x(t-h))\right)=&-hx^{(3)}(t)-\left(\frac{3h}{t}+1\right)x^{(2)}(t)\\
	&+\left(\frac{3h}{t^2}-\frac{3}{t}\right)x^{(1)}(t)+\mathcal{O}(h^2).
	\end{split}
	\end{align}
	We expand $x(t+h),\;x(t-h),\;x(t-2h)$ to the third order
	\begin{align*}
	\left(\alpha_1+\frac{\beta_1h}{t}+\frac{\gamma_1h^2}{t^2}\right)x(t+h)=&\left(\alpha_1+\frac{\beta_1h}{t}+\frac{\gamma_1h^2}{t^2}\right)\\&\left[x(t)+hx^{(1)}(t)+\frac{h^2}{2}x^{(2)}(t)+\frac{h^3}{6}x^{(3)}(t)+O(h^4)\right],\\
	\left(\alpha_3+\frac{\beta_3h}{t}+\frac{\gamma_3h^2}{t^2}\right)x(t-h)=&\left(\alpha_3+\frac{\beta_3h}{t}+\frac{\gamma_3h^2}{t^2}\right)\\&\left[x(t)-hx^{(1)}(t)+\frac{h^2}{2}x^{(2)}(t)-\frac{h^3}{6}x^{(3)}(t)+O(h^4)\right],\\
	\left(\alpha_4+\frac{\beta_4h}{t}+\frac{\gamma_4h^2}{t^2}\right)x(t-2h)=&\left(\alpha_4+\frac{\beta_4h}{t}+\frac{\gamma_4h^2}{t^2}\right)\\&\left[x(t)-2hx^{(1)}(t)+2h^2x^{(2)}(t)-\frac{4h^3}{3}x^{(3)}(t)+O(h^4)\right].
	\end{align*}
    Now we substitute these three equations and (\ref{fexp}) to truncation error of recurrence relation (\ref{ite2e})
	\begin{align*}
	\begin{split}
	L[x(t);h]=&\sum_{i=1}^4\left(\alpha_i+\frac{\beta_ih}{t}+\frac{\gamma_ih^2}{t^2}\right)x(t+(2-i)h)\\&+h^2\nabla F\left(x(t)+\frac{t-3h}{t}(x(t)-x(t-h))\right).
	\end{split}
	\end{align*} 
	Simple calculation shows that terms with order less than four will be eliminated if and only if the following equations are satisfied,
	\[\left\{\begin{aligned}\alpha_1&=2\\\alpha_2&=-5\\\alpha_3&=4\\\alpha_4&=-1\end{aligned}\right.,\qquad\left\{\begin{aligned}\beta_1&=\frac{9}{2}-k\\\beta_2&=-6+3k\\\beta_3&=\frac{3}{2}-3k\\\beta_4&=k\end{aligned}\right.,\qquad\left\{\begin{aligned}\gamma_1&=m_1\\\gamma_2&=-\frac{3m_1+m_2+3}{2}\\\gamma_3&=m_2\\\gamma_4&=\frac{m_1-m_2+3}{2}\end{aligned}\right.,\]
	where $k,\;m_1,\;m_2$ can be chosen randomly. Then we set $k=\frac{1}{2}$, $m_1=0$ and $m_2=3$ to get our SAG (\ref{ite2e}).
\end{proof}
\section{Algorithms for matrix completion experiment}\label{a5}
In this section we present FISTA (Algorithm \ref{a5F}), APG (Algorithm \ref{a5A}), FISTA with backtracking (Algorithm \ref{a5FB}) and APG with backtracking (Algorithm \ref{a5AB}).
\begin{algorithm}
	\caption{FISTA}
	\label{a5F}
	\begin{algorithmic}
		\STATE \textbf{Input: }step size $s$.
		\STATE \textbf{Initial value: }$\bm{Y}_1=\bm{X}_0=\bm{M}_{obs}$, $t_1=1$.
		\STATE \textbf{$\bm{k}$th iteration $(k\geq1)$.} Compute
		\begin{align*}
		&\bm{X}_k=\mathop{\arg\min}_{\bm{X}}\left\{\frac{1}{2s}\left\|\bm{X}-\left(\bm{Y}_k-sg(\bm{Y}_k)\right)\right\|^2+\lambda\|\bm{X}\|_*\right\},\\
		&t_{k+1}=\frac{1+\sqrt{1+4t_k^2}}{2},\\
		&\bm{Y}_{k+1}=\bm{X}_k+\frac{t_k-1}{t_{k+1}}(\bm{X}_k-\bm{X}_{k-1}).
		\end{align*}
	\end{algorithmic}
\end{algorithm}
\begin{algorithm}
	\caption{Accelerated proximal gradient method (APG)}
	\label{a5A}
	\begin{algorithmic}
		\STATE \textbf{Input: }step size $s$.
		\STATE \textbf{Initial value: }$\bm{X}_1=\bm{X}_0=\bm{M}_{obs}$.
		\STATE \textbf{$\bm{k}$th iteration $(k\geq1)$.} Compute
		\begin{align*}
		&\bm{Y}_k=\bm{X}_{k}+\frac{k-3}{k}(\bm{X}_{k}-\bm{X}_{k-1}),\\
		&\bm{X}_{k+1}=\mathop{\arg\min}_{\bm{X}}\left\{\frac{1}{2s}\left\|\bm{X}-\left(\bm{Y}_k-sg(\bm{Y}_k)\right)\right\|^2+\lambda\|\bm{X}\|_*\right\}.
		\end{align*}
	\end{algorithmic}
\end{algorithm}
\begin{algorithm}
	\caption{FISTA with backtracking}
	\label{a5FB}
	\begin{algorithmic}
		\STATE \textbf{Input: }some $\beta<1$.
		\STATE \textbf{Initial value: }$\bm{Y}_1=\bm{X}_0=\bm{M}_{obs}$, $t_1=1$, step size $s_0$.
		\STATE \textbf{$\bm{k}$th iteration $(k\geq1)$.}
		\STATE \hspace*{1em}Find the smallest nonnegative integer $i_k$ such that with $s=\beta^{i_k}s_{k-1}$
		\[F(\widetilde{\bm{X}})<F(\bm{Y}_k)+\left\langle\widetilde{\bm{X}}-\bm{Y}_k,g(\bm{Y}_k)\right\rangle+\frac{1}{2s}\|\widetilde{\bm{X}}-\bm{Y}_k\|^2,\]
		\hspace*{1em}where
		\[\widetilde{\bm{X}}=\mathop{\arg\min}_{\bm{X}}\left\{\frac{1}{2s}\left\|\bm{X}-\left(\bm{Y}_k-sg(\bm{Y}_k)\right)\right\|^2+\lambda\|\bm{X}\|_*\right\}.\]
		\STATE \hspace*{1em}Set $s_k=\beta^{i_{k}}s_{k-1}$ and compute
		\begin{align*}
		&\bm{X}_k=\widetilde{\bm{X}},\\
		&t_{k+1}=\frac{1+\sqrt{1+4t_k^2}}{2},\\
		&\bm{Y}_{k+1}=\bm{X}_k+\frac{t_k-1}{t_{k+1}}(\bm{X}_k-\bm{X}_{k-1}).
		\end{align*}
	\end{algorithmic}
\end{algorithm}
\begin{algorithm}
	\caption{APG with backtracking}
	\label{a5AB}
	\begin{algorithmic}
		\STATE \textbf{Input: }some $\beta<1$.
		\STATE \textbf{Initial value: }$\bm{X}_1=\bm{X}_0=\bm{M}_{obs}$, step size $s_1$.
		\STATE \textbf{$\bm{k}$th iteration $(k\geq1)$.} Compute
		\begin{align*}
		\bm{Y}_k=\bm{X}_{k}+\frac{k-3}{k}(\bm{X}_{k}-\bm{X}_{k-1}).
		\end{align*}
		\STATE \hspace*{1em}Find the smallest nonnegative integer $i_{k+1}$ such that with $s=\beta^{i_{k+1}}s_{k}$
		\[F(\widetilde{\bm{X}})<F(\bm{Y}_k)+\left\langle\widetilde{\bm{X}}-\bm{Y}_k,g(\bm{Y}_k)\right\rangle+\frac{1}{2s}\|\widetilde{\bm{X}}-\bm{Y}_k\|^2,\]
		\hspace*{1em}where
		\[\widetilde{\bm{X}}=\mathop{\arg\min}_{\bm{X}}\left\{\frac{1}{2s}\left\|\bm{X}-\left(\bm{Y}_k-sg(\bm{Y}_k)\right)\right\|^2+\lambda\|\bm{X}\|_*\right\}.\]
		\STATE \hspace*{1em}Set $s_{k+1}=\beta^{i_{k+1}}s_k$ and compute
		\begin{align*}
		\bm{X}_{k+1}=\widetilde{\bm{X}}.
		\end{align*}
	\end{algorithmic}
\end{algorithm}

\section{Detailed experiment results}\label{a6}
\begin{figure}[htbp]
	\centering
	\subfigure[$s=0.02$]{
		\includegraphics[width=4.3cm]{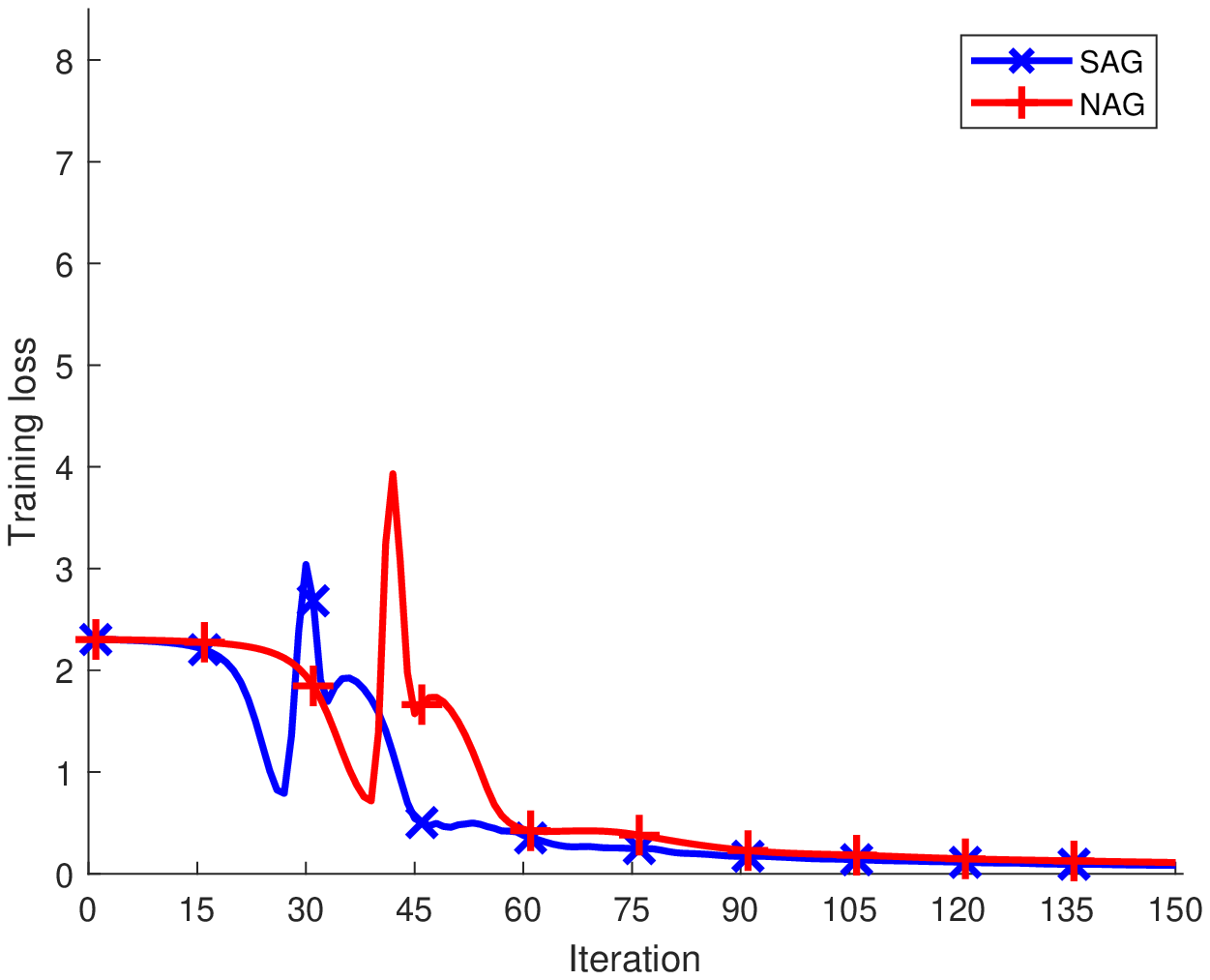}
	}
	\subfigure[$s=0.08$]{
		\includegraphics[width=4.3cm]{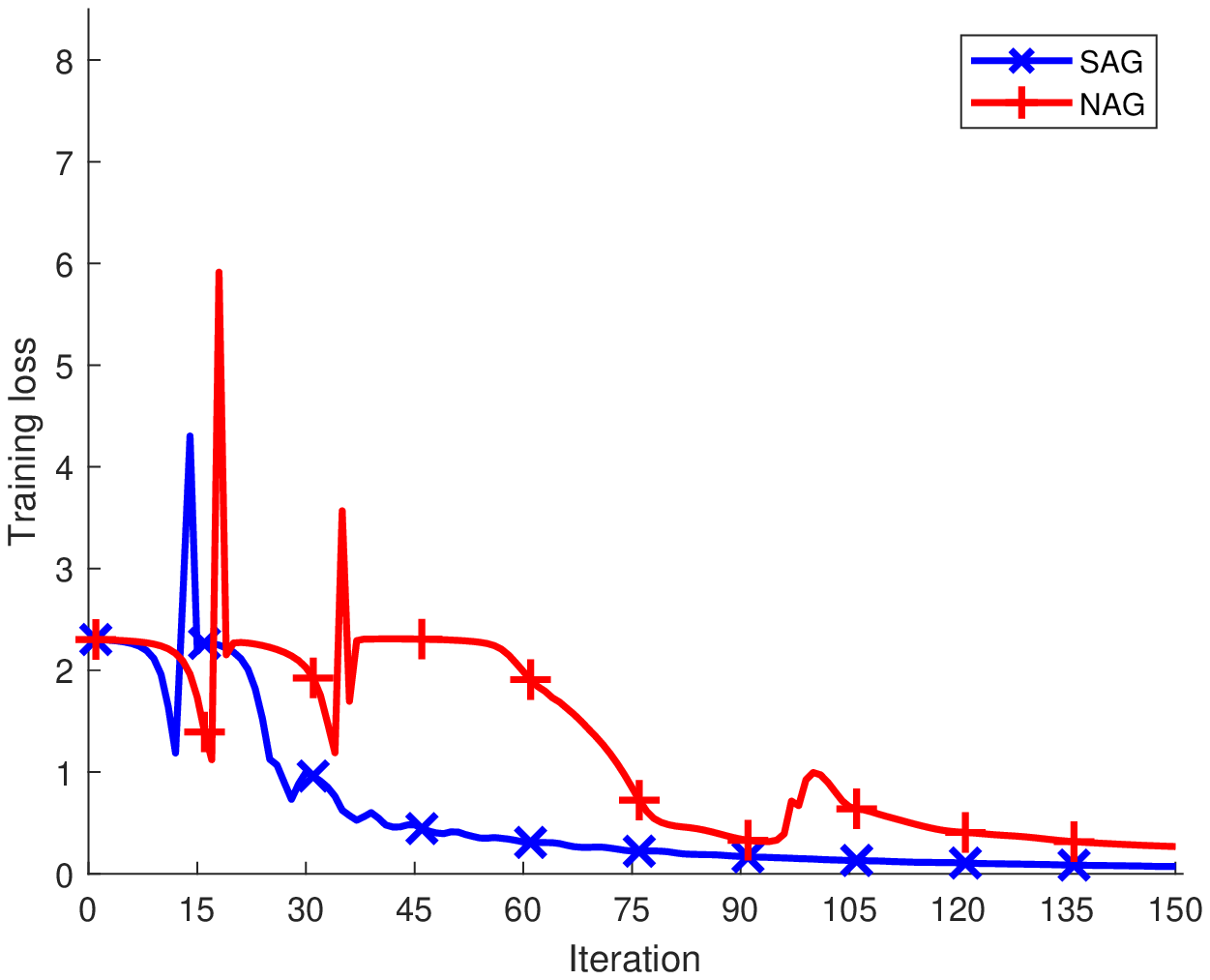}
	}
	\subfigure[$s=0.14$]{
		\includegraphics[width=4.3cm]{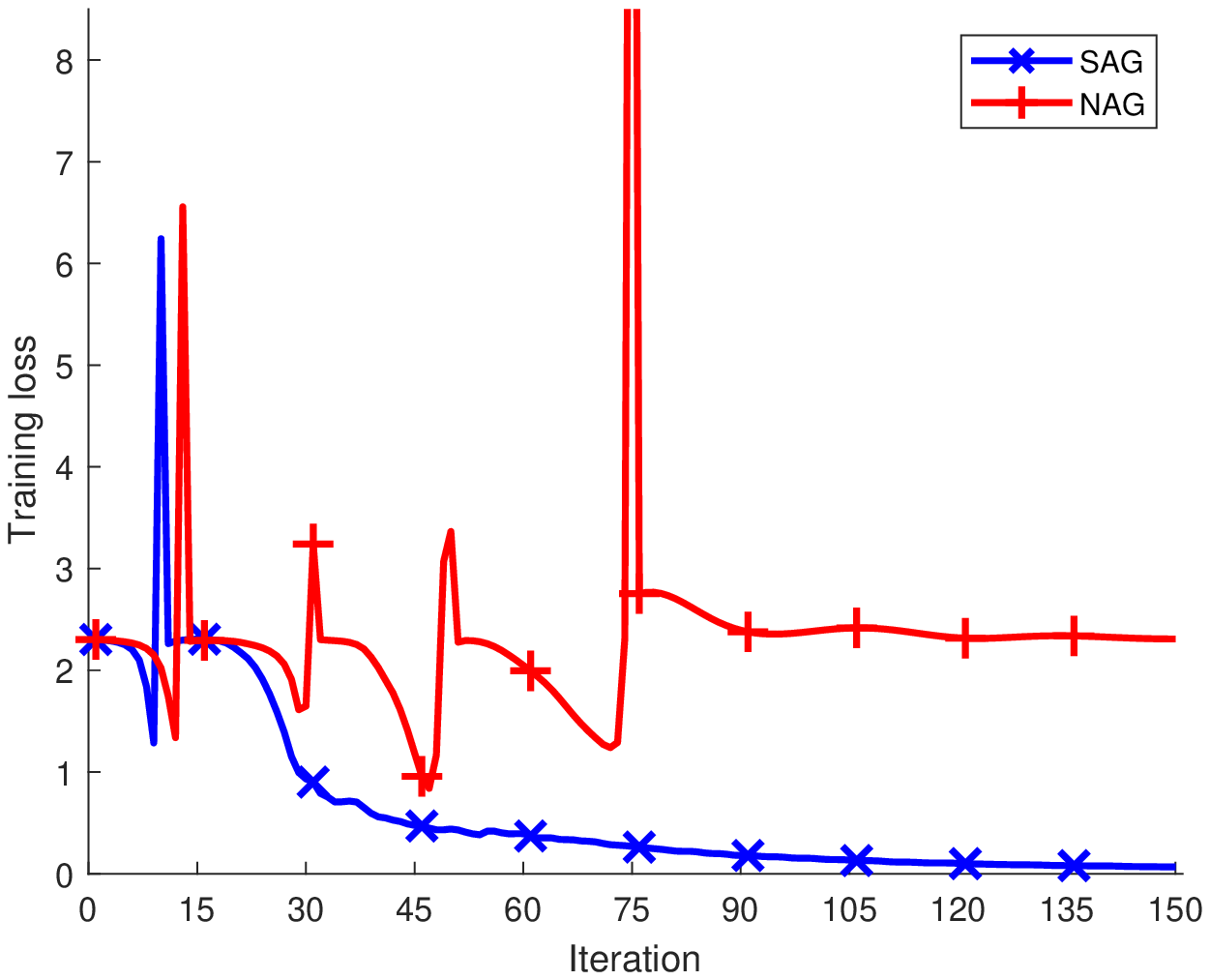}
	}
	\caption{Iterations of SAG and NAG for training loss. Y-axis represents the training loss for both methods. Step sizes for the three subgraphs are $0.02$, $0.08$ and $0.14$, coinciding with Figure \ref{fmni}.}
	\label{mlos}
\end{figure}

We use cross-entropy loss to train the network for handwriting digit recognition task. In every iteration of the two methods, we calculate the gradient using all the $20000$ training images, so it is meaningful to study the change of the training loss during iteration. In Figure \ref{mlos} we plot the training loss against the number of iterations for SAG and NAG. This figure also shows that SAG is more stable than NAG when step size is large.
\end{document}